\documentclass[a4paper,leqno]{article} 
\pdfoutput=1
\usepackage[utf8]{inputenc}
\usepackage[T1]{fontenc}
\usepackage{enumitem,amsfonts,amsmath,amssymb,amsthm,stmaryrd}
\usepackage{lmodern} 
\usepackage{dsfont} 
\usepackage{fullpage} 
\usepackage{microtype} 
\usepackage{comment} 
\usepackage{graphicx}
\usepackage{mathrsfs}
\usepackage{tikz} 
\usepackage{template}
\usepackage{float} 
\usepackage{authblk}
\DeclareMathOperator\sP{P}
\DeclareMathOperator\E{E}

\title{Law of large numbers for a finite-range random walk in a dynamic random environment with nonuniform mixing}
\author{Julien \textsc{Allasia}\thanks{Universite Claude Bernard Lyon 1, CNRS, Ecole Centrale de Lyon, INSA Lyon, Université Jean Monnet, ICJ UMR5208,
69622 Villeurbanne, France. \textit{\href{mailto:allasia@math.univ-lyon1.fr}{allasia@math.univ-lyon1.fr}}.}}
\date{}

\begin{document}
\maketitle

\begin{abstract}
In this paper, we study random walks evolving on $\ZZ$ in a dynamic random environment that we assume to have time correlations that decrease polynomially fast. We show a law of large numbers by generalizing methods already used for the nearest-neighbor framework to the finite-range one. This requires some new ideas to get around the absence of a monotonicity property that was crucial in the proof for the nearest-neighbour case. Our proof works both in discrete and continuous time.
\end{abstract}

\begin{abstract}
Dans ce papier, nous étudions des marches aléatoires qui évoluent sur $\ZZ$ dans un milieu aléatoire dynamique dont les corrélations en temps sont supposées décroître de façon polynomiale. On montre une loi des grands nombres en généralisant des méthodes déjà utilisées pour le modèle aux plus proches voisins au modèle à portée finie. Cela requiert de nouvelles idées pour contourner l'absence d'une propriété de monotonie qui était cruciale dans la preuve dans le modèle aux plus proches voisins. Notre preuve fonctionne à la fois en temps discret et en temps continu.
\end{abstract}

\section{Introduction}

The study of random walks in random environments (RWRE) has started in the 60s with a one-dimensional static model aimed at studying DNA replication in \cite{Chernov}. The model was the following: we let $(\alpha_i)_{i\in\ZZ}$ be a sequence of random variables in $[0,1]$. Conditioned on $(\alpha_i)$, we define a random walk on environment $(\alpha_i)$ as a Markov process $(X_n)_{n\in\NN}$ starting at $0$ with transition matrices $M(i,i+1)=\alpha_i$ and $M(i,i-1)=1-\alpha_i$. Later, this model found further applications in metallurgy (\cite{Temkin}), biology (\cite{CoccoM}) and physics (\cite{LDMF}, \cite{Hughes}).

The law of large numbers (LLN) for this model under ergodic assumptions on the environment was proven in \cite{Solomon}, and today this model is thoroughly understood (see for instance \cite{KKS}, \cite{Sinai}, \cite{Shi} and \cite{ESZ} for further results). If we lift the nearest-neighbor assumption to consider finite-range RWRE in dimension 1, the law of large numbers also holds, as is shown in \cite{Bremont1} and \cite{Bremont2}.

However, generalizing to higher dimensions or dynamic environments (the environment is allowed to evolve in time) still requires some strong assumptions. For a review on the static setting, we refer to \cite{Sznitman_notes}, \cite{Zeitouni} and \cite{Revesz}.

In the dynamic framework, we allow the environment to change at each time step. In order to get asymptotic results, strong assumptions on this dynamics are often made, even in dimension $1$, starting with independence (see \cite{BoldMinPel}, \cite{RAS}). Weaker mixing conditions have also been investigated, but usually the mixing has to be either very fast or uniform (see \cite{CZ}, \cite{KO}, \cite{AHR}, \cite{ABF}). Progress has been made for environments which do not satisfy these strong mixing conditions, but the techniques used are often highly dependent on the specificities of the environments considered. Examples of such environments include the contact process (\cite{HollSan}, \cite{MV}, \cite{Bet}), the exclusion process (\cite{ASV}, \cite{HuvSim}) and Poissonian fields of random walks (\cite{HKS}, \cite{HHSST}, \cite{BHSSST1}, \cite{BHSSST2}).

In \cite{BHT} however, the authors manage to adapt techniques of multi-scale renormalization to show a LLN in the general framework of one-dimensional environments having time correlations that decrease polynomially fast. This framework encompasses a lot of models that were studied before, however the techniques used here only yield a LLN.

In this paper, our goal is to generalize the methods used to show the LLN in \cite{BHT} to the finite-range framework, in which random walks are allowed to jump to all the sites in $\ZZ$ at distance at most $R\geqslant 1$. Under assumptions of translation invariance and mixing of the environment as well as uniform ellipticity of the random walk, we show there exists $\nu\in [-R,R]$ such that almost surely, $X_t/t$ goes to $\nu$ as $t\rightarrow\infty$.

Our methods work for both discrete and continuous time, therefore they also allow for a generalization of the LLN in \cite{BHT} to the discrete-time framework, which is not straightforward using the methods of \cite{BHT}.

The problem with the methods from \cite{BHT} in our finite-range framework is that one of the key arguments in \cite{BHT} is a monotonicity property in the coupling of random walks that precisely relied on the nearest-neighbor assumption (and, to a lesser extent, the continuous-time framework). In order to adapt the ideas of traps from \cite{BHT}, we add a classic uniform ellipticity assumption that will allow us to control the probability for two paths to coalesce, which will help us recover an essential part, albeit weaker than monotonicity, of the behavior of the coupling in \cite{BHT}.

This generalization of the innovative techniques used in \cite{BHT} shows that they are not as rigid as they may have seemed at first glance and opens up the possibility of other generalizations. In \cite{article_stat} for example, we move from the dynamic one-dimensional framework to a static two-dimensional one and we show the existence of an asymptotic direction for a ballistic random walk.

\paragraph{Outline of the paper.}
In Section \ref{D:s:framework}, we define precisely the framework of this paper by defining one-dimensional dynamic environments and random walks on them, before stating our main results, Theorems \ref{D:t:LLNd} and \ref{D:t:LLNc}. Its proof is divided into two parts, which correspond respectively to Sections \ref{D:s:limiting_speeds} and \ref{D:s:v_-=v_+}. In the first part, we show the existence of two limiting speeds that bound the average speed of our random walk over long time intervals. In the second part, we show that these two speeds coincide, which will give the limiting speed that we are after for the LLN. The final proof of our theorems, using results shown later in the paper, can be found in Section \ref{D:ss:final_proof}, both in discrete and continuous time.

\paragraph{Conventions.}\hfill

$\N$, $\Z$ and $\R$ respectively denote the set of natural integers (starting from $0$), relative integers and real numbers. $\N^*$ denotes $\N\setminus \{0\}$, $\RR_+$ denotes $\{x\in\RR,\,x\geqslant 0\},$ $\RR_+^*$ is $\RR_+\setminus\{0\}.$ If $n\leqslant m$ are two integers, $\llbracket n,m\rrbracket$ is the set of integers $[n,m]\cap\Z.$ For a couple $x=(a,b)\in\RR^2$, we write $a=\pi_1(x)$ and $b=\pi_2(x)$, and we refer to them as the horizontal and vertical coordinates of $x$. If $S$ is a finite set, $|S|$ and $\# S$ denote the cardinality of $S$.

$c$ denotes a positive constant that can change throughout the paper and even throughout one single computation. Constants that are used again later in the paper will be denoted with an index when they appear for the first time.

Drawings across the paper are not to scale and they do not necessarily represent the random walks in a perfectly accurate way: they are only meant to make the reading easier. For instance, sample paths will often look like smooth curves, although our random walks evolve on $\ZZ.$

\section{Framework}\label{D:s:framework}

\subsection{Setup of the problem and main result}
Let $\mathbb{T}$ be either $\NN$ or $\RR_+$ depending on if we work in discrete or continuous time. Let $\LLL=\ZZ\times\mathbb{T}$, the set of space-time points. We want to study the asymptotic behavior of a random walk moving randomly in a dynamic random environment on $\ZZ$. More precisely, we want to have an asymptotic speed for this random walk. In order to do this, it is very useful to couple random walks starting at different space-time points in $\LLL.$

Let $S$ be a topological space (in common applications, $S$ is usually a finite or countable subset of $\RR$). Let $(\Omega,\mathcal{T},\sP)$ be a probability space. A random environment defined on {\color{black}$(\Omega,\mathcal{T},\sP)$} is a random variable taking values in the set of càdlàg functions from $\mathbb{T}$ to $S^\ZZ$, endowed with the product topology on $S^\ZZ$ and the subsequent Borel $\sigma$-algebra (note that in discrete time, the càdlàg assumption is useless). For $t\in\mathbb{T}$ and $x\in\ZZ$, $\omega_t(x)$ is called the state of environment $\omega$ at time $t$ and site $x$.

We now assume that environment $\omega$ is given and we give a formal construction of the random walk in $\omega$. What we want to do is to couple random walks starting from different points in $\LLL$ together. In order to do so, we use uniform random variables to encapsulate the randomness needed for the jumps of the random walks and we allocate them to points in $\LLL$.

Let $\ell,R\in\N^*$. We give ourselves a {\color{black}jump} function, that is a measurable function $g:S^{\llbracket -\ell,\ell\rrbracket}\times [0,1]\to \llbracket -R,R\rrbracket$. We call $\ell$ the dependency range and $R$ the {\color{black}jump} range of the random walk. We also define a family of independent uniform random variables in $[0,1]$, defined in a probability space $(\Omega',\mathcal{T}',\sP')$, denoted by $(U_n^x)_{n\in\N,x\in\Z}$. Those will give the extra randomness that we need for the jumps. We set
\begin{align*}
\Omega_0=\Omega\times\Omega',\;\; \mathcal{T}_0=\mathcal{T}\otimes\mathcal{T}',\;\; \PP=\sP\otimes \sP'.
\end{align*}
\begin{definition}[discrete-time setting]\label{D:d:random_walk_discrete}
Suppose $\mathbb{T}=\NN$. Let $y=(x_0,n_0)\in \LLL$. The random walk in environment $\omega$ starting at $y$ is the random process on $(\Omega_0,\mathcal{T}_0)$ given by $n\in\NN\mapsto Z_n^y=(X_n^y,n_0+n)\in\LLL$ such that $$\left\{\begin{array}{l}X_0^y=x_0;\\
    \forall n\in\NN,\; X^y_{n+1}=X^y_n+g\left(\omega_{n_0+n}(X^y_n-\ell),\ldots,\omega_{n_0+n}(X^y_n+\ell),U_{n_0+n}^{X_n^y}\right).\end{array}\right.$$
\end{definition}

For the continuous-time setting, we simply define the random walk starting at the origin $(0,0)$. Indeed, defining $(X_n^y)$ for any starting point $y$ in Definition \ref{D:d:random_walk_discrete} was only meant to be used for the proof of the LLN, and the LLN in continuous time will be a consequence of the LLN in discrete time.

\begin{definition}[continuous-time setting]\label{D:d:random_walk_continuous}
Suppose $\mathbb{T}=\RR_+$. Let $(T_n)_{n\in\NN^*}$ be a Poisson process in $\RR_+^*$ of parameter $1$, independent of $(\omega,(U_n^x)_{n\in\NN,x\in\ZZ})$. The random walk in environment $\omega$ starting at the origin is the {\color{black}càdlàg} random process on $(\Omega_0,\mathcal{T}_0)$ given by $t\in\RR_+\mapsto Z_t=(X_t,t)\in\LLL$ such that $$\left\{\begin{array}{l}X_0=0;\\
    \forall t>0,\;X_t=
    \left\vert\begin{array}{ll}
    X_{t^-}&\mbox{if $t\notin\{T_n,n\in\NN^*\}$,}\\
    X_{t^-}+g\left(\omega_{T_n}(X_{t^-}-\ell),\ldots,\omega_{T_n}(X_{t^-}+\ell),U_{n-1}^{X_{t^-}}\right)&\mbox{if $t=T_n.$}\end{array}\right.\end{array}\right.$$
\end{definition}

The first important thing to say about those definitions is that they define random walks in random environments in the usual sense, namely, when the environment is fixed, $X$ is a Markov process. Furthermore, the coupling given by Definition \ref{D:d:random_walk_discrete} implies a crucial property for the random walks. The definition of $X^y$ uses the same environment $\omega$ and the same variables $(U_n^x)$ for the jumps, independently of $y$. This implies that two random walks that meet at the same space-time point will then have the exact same sample paths. Here it is crucial that we allocated the randomness needed for the jumps (the $(U_n^x)$ variables) to space-time points instead of associating one sequence of uniform variables to each random walk. More precisely, we have the following coalescing property:
\begin{align}\label{D:e:coupling_property}
    \forall y,y'\in\LLL,\; \forall n,n'\in\NN,\;Z^y_n=Z^{y'}_{n'}\Longrightarrow \forall k\in\NN,\,Z^y_{n+k}=Z^{y'}_{n'+k}.
\end{align}
\begin{remark}\label{D:r:monotonicity_BHT}
This property has a much stronger version in \cite{BHT}. Under the assumptions therein, particularly the fact that $R=1$, sample paths cannot cross without meeting {\color{black}(meaning that if the two segmented lines in $\R\times\R_+$ obtained by joining the points in $\LLL$ visited by the sample paths cross, then the sample paths visit
some common point in $\LLL$)}. This implies, using \eqref{D:e:coupling_property}, that if a random walk starts on the left of another random walk, it cannot end up on its right. This monotonicity property (see (2.9) in \cite{BHT}) was crucial in the proof of the LLN in \cite{BHT}, and it is lost when $R\geqslant 2.$ We refer to Section \ref{D:ss:comparison} for a detailed discussion about the differences between our paper and \cite{BHT}.
\end{remark}

For the rest of the paper, we make the following assumptions on the environment.
\begin{assumption}[Translation invariance]\label{D:a:invariance}
For every $(x_0,t_0)\in\LLL$, 
\begin{align*}
(\omega_t(x))_{(x,t)\in\LLL}\text{ and }(\omega_{t_0+t}(x_0+x))_{(x,t)\in\LLL}\text{ have the same law.}
\end{align*}
\end{assumption}

For the next assumption, we need some definitions. Let $B\subseteq\RR^{2}$. We say that $B$ is a box if it is of the form $[a_1,b_1)\times[a_2,b_2)$ with $a_i<b_i$. If $B$ is a box, we define its spatial diameter $\mathrm{diam}(B)=b_1-a_1$ and height $h(B)=b_2-a_2$. If $B=[a_1,b_1)\times[a_2,b_2)$ and $B'=[a_1',b_1')\times[a_2',b_2')$ are two boxes, we call $\mathrm{sep}(B,B')$ their vertical separation defined by
\begin{align*}
 \mathrm{sep}(B,B')= \left\{\begin{array}{ll}
              a_{2}-b'_{2}& \text{ if } b'_2< a_2,\\
              a'_2-b_2&\text{ if }b_2<a'_2,\\
              0&\text{ else.}
             \end{array}\right.
\end{align*}
{\color{black}
We also set, for $A\subset \LLL$, $\omega\vert_A=\{\omega_t(x),(x,t)\in A\}$.}

\begin{definition}\label{D:d:mixing}
\ncD{D:c:mixing} Let $\ucD{D:c:mixing},a,b,\alpha>0.$ We denote by $\mathcal{D}_{\mathbb{T}}(\ucD{D:c:mixing},a,b,\alpha)$  the following property: for every $H>0,$ for every pair of boxes $B,B'\subseteq\RR^2$ satisfying
$$\max(\mathrm{diam}(B),\mathrm{diam}(B'))\leqslant aH,\quad \max(h(B), h(B'))\leqslant bH,\quad \mathrm{sep}(B,B')\geqslant H,$$
and for every pair of $\{0,1\}$-valued functions $f_1$ and $f_2$ on $(S^\ZZ)^{\mathbb{T}}$ such that $f_1(\omega)$ is $\sigma(\omega\vert_{B\cap\LLL})$-measurable and $f_2(\omega)$ is $\sigma(\omega\vert_{B'\cap\LLL})$-measurable,
    \begin{align*} \E[f_1(\omega)f_2(\omega)]\leqslant \E[f_1(\omega)]\E[f_2(\omega)]+\ucD{D:c:mixing}\, H^{-\alpha}.
    \end{align*}
\end{definition}

Note that we did not require a lower bound for the covariance. The last assumption we make is a standard uniform ellipticity assumption that we formulate using the {\color{black}jump} function.
\begin{assumption}[Uniform ellipticity]\label{D:a:uniform_ellipticity}
There exists a constant $\gamma>0$ such that if $U$ is uniformly distributed in $[0,1]$ under $\PP$,
\begin{align*} \inf_{x\in\llbracket -R,R\rrbracket}\inf_{\sigma_{-\ell},\ldots,\sigma_\ell\in S} \PP(g(\sigma_{-\ell},\ldots,\sigma_\ell,U)=x)\geqslant \gamma.
\end{align*}
\end{assumption}

When the starting point $y$ of random walk $X^y$ is omitted, we mean $y=(0,0)$. The goal of this paper is to show the following two theorems, the first one dealing with the discrete-time case ($\mathbb{T}=\NN$) and the second one with the continuous-time one ($\mathbb{T}=\RR_+$).

\ncD{D:c:LLN}
\begin{theorem}[LLN in discrete time]\label{D:t:LLNd}
Assume that Assumptions \ref{D:a:invariance} and \ref{D:a:uniform_ellipticity} are satisfied, as well as $\mathcal{D}_{\NN}(\ucD{D:c:mixing},2R+3,1,\alpha)$ for some parameters $\ucD{D:c:mixing}>0$ and $\alpha>11$. Then there exists $\nu\in [-R,R]$ such that \begin{align*}
    \PP\text{-almost surely},\;\;\frac{X_n}{n}\xrightarrow[n\to\infty]{}\nu.
\end{align*}
Moreover we have a polynomial rate of convergence:
\begin{align}\label{D:LLN_proba}
\forall \varepsilon>0,\;\exists\, \ucD{D:c:LLN}=\ucD{D:c:LLN}(\varepsilon)>0,\;\forall n\in\NN^*,\;\;\;\PP\left(\left\vert \frac{X_n}{n}-\nu\right\vert\geqslant\varepsilon\right)\leqslant \ucD{D:c:LLN}\,n^{-\alpha/4}.\end{align}
\end{theorem}

\begin{theorem}[LLN in continuous time]\label{D:t:LLNc}
Assume that Assumptions \ref{D:a:invariance} and \ref{D:a:uniform_ellipticity} are satisfied, as well as $\mathcal{D}_{\RR_+}(\ucD{D:c:mixing},a,b,\alpha)$ for some parameters $\ucD{D:c:mixing}>0$, $a>2R+3$, $b>1$ and $\alpha>11$. Then there exists $\nu\in [-R,R]$ such that \begin{align*}
    \PP\text{-almost surely},\;\;\frac{X_t}{t}\xrightarrow[t\to\infty]{} \nu.
\end{align*}
Moreover we have a polynomial rate of convergence:
\begin{align*}
\forall \varepsilon>0,\;\exists\, \ucD{D:c:LLN}=\ucD{D:c:LLN}(\varepsilon)>0,\;\forall t>0,\;\;\;\PP\left(\left\vert \frac{X_t}{t}-\nu\right\vert\geqslant\varepsilon\right)\leqslant \ucD{D:c:LLN}\,t^{-\alpha/4}.\end{align*}
\end{theorem}

Note that the lower bound on $\alpha$ may not be optimal. However, in practice most examples of environments that we use satisfy $\mathcal{D}_{\mathbb{T}}(\ucD{D:c:mixing},a,b,\alpha)$ with arbitrarily high $a$, $b$ and $\alpha$, for a proper choice of $\ucD{D:c:mixing}$, so we are not too interested in optimality here.

\subsection{Examples}\label{D:ss:examples}
The strength of the methods developed here is that they apply to a large class of environments. Examples of such environments can be found in \cite[page 2041]{BHT}. Mind that they are defined in a continuous-time setting, but we can get discrete-time environments simply by considering their values on times in $\NN$. More specifically, whenever a continuous-time environment $(\eta_t(x))_{t\in\RR_+,x\in\ZZ}$ from \cite{BHT} satisfies the property called $\mathcal{D}_{\text{env}}(c_0,\alpha)$ for some $\alpha$, then the discrete-time environment defined by $\omega_n(x)=\eta_n(x)$ for $n\in\NN$ and $x\in\ZZ$ satisfies our property $\mathcal{D}_\NN(c_0,5,1,\alpha)$, and taking a closer look at the proofs of $\mathcal{D}_{\text{env}}(c_0,\alpha)$ for those examples, we see that there is no obstacle to getting $\mathcal{D}_\NN(c_0,2R+3,1,\alpha)$ with $\alpha>11.$ The examples of such environments given in \cite{BHT} are the contact process, independent renewal chains and Markov processes with a positive spectral gap (for instance, the East model).

\subsection{Comparison with \cite{BHT}}\label{D:ss:comparison}

Let us compare our results to those of \cite{BHT} in order to better understand the specificities of our paper. There are two main differences in the assumptions made. One involves where the random walk can jump and the other one the times at which the random walk jumps. 

The first one is linked to the core interest of our paper, namely generalizing the methods of \cite{BHT} to prove a LLN when the random walk is finite-range. In \cite{BHT}, the random walk is assumed to satisfy $R=1$ (note that the random walk is allowed to remain at the same site at a {\color{black}jump} time). Here however, we allow our random walk to jump in $\llbracket -R,R\rrbracket$ for $R\geqslant 2$. {\color{black}Nonetheless, as a compensation, we add a uniform ellipticity assumption (cf. Assumption \ref{D:a:uniform_ellipticity}) that was absent in \cite{BHT}.}

The second difference is trickier: it involves the jump times of the RWRE. First of all, \cite{BHT} only deals with a continuous-time framework (but we still denote by $\omega$ the environment). This is somewhat restrictive, but the assumptions on the {\color{black}jump} times {\color{black}made in \cite{BHT}} are quite weak. {\color{black}Let us recall those assumptions. For each $x\in\ZZ$, we define the sequence of jump times at site $x$ to be an increasing sequence $(T_i^x)_{i\in\NN^*}$ of random times in $\RR_+$} such that almost surely, for every $x\in\ZZ$,
\begin{align}\label{D:a:ass_times}
    \{T_i^x,\,i\in\NN^*\}\text{ and }\{T_i^{x+1},\,i\in\NN^*\}\text{ are disjoint.}
\end{align}

Moreover, some control on the number of {\color{black}jump} times that fall in a certain time interval is also assumed, in order for the random walks not to jump unreasonably often in this time interval. Without such control, proving a LLN is of course hopeless. A tool to formulate this kind of assumption is the notion of allowed path. An allowed path on $[0,T]$ is a function $t \in [0,T]\mapsto \beta(t)=(\beta_1(t),\beta_2(t))\in\LLL$ (here, $\LLL=\ZZ\times\RR_+$ just as in the continuous-time framework presented in Definition \ref{D:d:random_walk_continuous}), càdlàg, with jumps in $\{-1,0,1\}$, such that
$$\left\{\begin{array}{l}
\forall s\geqslant 0,\;\;\beta_2(t+s)-\beta_2(t)=s;\\
\text{if $\beta(t)=(x,s)$, then $\beta_1(t+r)=x$ for every ${\color{black}r<\min\{T_i^x-s,i\text{ such that }T_i^x>s\}.}$}\end{array}\right.$$
{\color{black} Here, the assumption made in \cite{BHT}} is a control on the speed that allowed paths can attain and their localization in a box until a certain time: for every $v>1,$
\ncD{D:c:BHTbox}
\begin{align}
    &\liminf_{T\to\infty} \PP\left(\begin{array}{c}\text{$\exists\,\beta$ allowed path on $[0,T]$, starting in $[0,T)\times\{0\}$}\\ \text{and such that $\beta_1(T)-\beta_1(0)\geqslant vT$}\end{array}\right)=0\;\\
    &\liminf_{T\to\infty} \PP\left(\begin{array}{c}\text{$\exists\,\beta$ allowed path on $[0,T]$, starting in $[0,T)\times\{0\}$}\\ \text{and such that $\beta_1(T)-\beta_1(0)\leqslant -vT$}\end{array}\right)=0,\end{align}
and there exists $\ucD{D:c:BHTbox}>0$ such that for all $T>0$,
\begin{align}
\PP\left(\begin{array}{c} \text{$\exists\,\beta$ allowed path on $[0,T]$, starting at $0$,} \\ \text{$\exists\, s\in[0,T],\;[\beta_1(s)-\ell,\beta_1(s)+\ell] \nsubseteq [-2T,2T]$} \end{array}\right)\leqslant \ucD{D:c:BHTbox}^{-1}\,e^{-\ucD{D:c:BHTbox}T}.\label{large_deviation}\end{align}

\ncD{D:c:BHTdec}
{\color{black}
Now, the main assumptions on the environments required to show asymptotic results with the tools introduced in \cite{BHT} are translation-invariance and polynomial mixing. Since the assumptions on the {\color{black}jump} times recalled so far are quite weak, they need to include the jump times in the invariance and mixing assumptions: they assume that the process given by $(\omega,\mathbf{T})$ where $\mathbf{T}=(T_i^x)_{x\in\ZZ,i\in\NN^*}$ is translation-invariant and that it mixes polynomially. Let us formulate this rigorously. First, if we set, for $(z,s)\in\LLL$, $j(z,s,x)=\min\{i\in\N^*,\,T_i^{z+x}>s\}$, we ask that
\begin{align}\forall (z,s)\in\LLL,\,(\omega,\mathbf{T})\text{ and }\left((\omega_{s+t}(z+x))_{(x,t)\in\LLL},\,(T_{j(z,s,x)+i-1}^{z+x}-s)_{x\in\ZZ,\,i\in\N^*}\right)\text{ have same law,}\end{align}
Then, for $A\subset\LLL$, we define $\sigma_A$ to be the sigma-algebra generated by $\omega\vert_A$ and the times $(T_i^x)_{x,i}$ such that $(x,T_i^x)\in A.$ The assumption is that there exists $\ucD{D:c:BHTdec}
>0$ and $\alpha>8$ such that for every pair of boxes $B,B'\subseteq\RR^2$ satisfying $$\max(\mathrm{diam}(B),\mathrm{diam}(B'))\leqslant 5H,\quad \max(h(B), h(B'))\leqslant H,\quad \mathrm{sep}(B,B')\geqslant H,$$
and for every pair of $\{0,1\}$-valued functions $f_1$ and $f_2$ on $(S^\ZZ)^{\mathbb{\RR_+}}\times \RR_+^{\ZZ\times\NN^*}$ such that $f_1(\omega,\mathbf{T})$ is $\sigma_B$-measurable and $f_2(\omega,\mathbf{T})$ is $\sigma_{B'}$-measurable,
\begin{align}\label{D:e:mix_BHT}
\EE[f_1(\omega,\mathbf{T})\,f_2(\omega,\mathbf{T})]\leqslant\EE[f_1(\omega,\mathbf{T})]\,\EE[f_2(\omega,\mathbf{T})]+\ucD{D:c:BHTdec}\,H^{-\alpha}.
\end{align}}

Note that assumptions \eqref{D:a:ass_times} to \eqref{D:e:mix_BHT} are satisfied if the $((T_i^x)_{i\in\NN^*},x\in\ZZ)$ are independent Poisson processes of parameter $1$, drawn independently of $\omega$. This is very similar to what we assume in our continuous-time setting: see Definition \ref{D:d:random_walk_continuous}. However, in this paper, we make a strong use of the Poisson behavior and of the independence between the {\color{black}jump} times and the environment, so we cannot weaken the assumptions on the {\color{black}jump} times so much: see Section \ref{D:sss:final_proof_continuous}.

The proof of Theorems \ref{D:t:LLNd} and \ref{D:t:LLNc} is more subtle than for the LLN in \cite{BHT}. Indeed, the fact that random walks can jump at larger range allows them to cross paths without meeting. On the contrary, in \cite{BHT}, the $R=1$ assumption, along with \eqref{D:a:ass_times}, forces the random walks to keep their initial order: if a random walk starts on the left of another random walk, it will never end up on its right - in the worst case scenario, they will coalesce (see Remark \ref{D:r:monotonicity_BHT} in \cite{BHT}). Our goal will be to recover that behavior by forcing random walks to coalesce when they become dangerously close to each other. The need for a lower bound on the probability of this scenario of coalescence is the reason why we made a uniform ellipticity assumption. Since the monotonicity argument from \cite{BHT} falls apart in any case, \eqref{D:a:ass_times} is not needed anymore, and this is why our proof actually works in discrete time without any problem.

Actually, although it is not mentioned in \cite{BHT}, the methods used therein can be adapted for the discrete-time setting, {\color{black}provided that the {\color{black}jump} function $g$ takes values in $\{-1,1\}$ (we call this a nearest-neighbor setting)}. The trick is to change $\LLL=\ZZ\times\NN$ to $\LLL=\{(x,n)\in\ZZ\times\NN,\,x+n=0\,[2]\}$, and only consider random walks starting in $\LLL$. Those have their sample paths in $\LLL$ by construction. Now, two allowed sample paths on $\LLL$ cannot cross without meeting, so we recover the monotonic behavior.

The following table provides a recap of the above discussion, by displaying, for each model, the assumptions on the {\color{black}jump} times, on the jumps of the random walk, the methods used to show the LLN and the set $\LLL$ on which the random walks evolve.
{\small
\begin{center}
\begin{tabular}{|c|c|c|c|c|} 
\hline
Model & \cite{BHT} & Discrete adaptation & Continuous finite-range & Discrete finite-range\\
\hline
Time & \eqref{D:a:ass_times} to \eqref{D:e:mix_BHT} & $T_i^x=i$ & $\mathcal{P}(1)$ i.i.d. independent of $\omega$& $T_i^x=i$ \\
\hline 
Jumps & In $\{-1,0,1\}$ & In $\{-1,1\}$ & \multicolumn{2}{|c|}{In $\llbracket -R,R\rrbracket$, uniform ellipticity}\\
\hline
Methods used & \multicolumn{2}{|c|}{\cite{BHT}} & \multicolumn{2}{|c|}{Present paper}\\
\hline
Set $\LLL$ & $\ZZ\times\RR_+$ & $\{(x,n),\,x+n=0[2]\}$ & $\ZZ\times\RR_+$ & $\ZZ\times\NN$\\
\hline\end{tabular}\end{center}}

As far as the proof of the LLN is concerned, we refer the reader to Section \ref{sss:roadmap} for a better understanding of the differences between ours and that of \cite{BHT}.

To finish our comparison, let us give a warning about an interesting application of the tools developed in \cite{BHT} that is not covered by this paper (see \cite{BHT}, Section 8.3). Recall that in \cite{BHT}, no assumption of independence between $\omega$ and the $(T_i^x)$ is made. Therefore, the {\color{black}jump} times might actually be constructed using times involved in the construction of the environment itself. This is what is done with the East model in order to show Proposition 8.5 (page 2044), which describes the asymptotic behavior of the distinguished zero, viewed as a random walk of range $1$ on the East model. The random times used in the graphical construction of the latter are exactly the {\color{black}jump} times for this RWRE. Therefore, this result is not a consequence of the results shown in this paper. Actually, there is a second reason for this: the distinguished zero does not satisfy our uniform ellipticity assumption (actually, the distinguished zero can never go to the left by construction). That being said, the East model is a random environment to which our results can apply, provided that the {\color{black}jump} times chosen for the random walk satisfy the assumptions of Definition \ref{D:d:random_walk_continuous}, since it is a Markov process with positive spectral gap: see Section \ref{D:ss:examples}.

The rest of the paper is dedicated to proving Theorem \ref{D:t:LLNd} {\color{black}(in Section \ref{D:sss:final_proof_continuous} we show that Theorem \ref{D:t:LLNc} can be derived from Theorem \ref{D:t:LLNd}).} Its final proof using results that will be shown in the next two sections can be found in Section \ref{D:ss:final_proof}. Mind that except in Section \ref{D:sss:final_proof_continuous}, where we will deduce the LLN in continuous time from that in discrete time, we only work with discrete time: $\mathbb{T}=\NN$ and $\LLL=\ZZ\times\NN$.

\subsection{Tools of the proof}\label{D:ss:strategy}

\subsubsection{Mixing property}\label{D:ss:dec}
Figure \ref{D:f:mixing} below provides an illustration for assumption $\mathcal{D}_{\NN}(\ucD{D:c:mixing},2R+3,1,\alpha)$. In practice, we will use mixing for events that involve our random walks and not only our environment. This is actually not stronger than the mixing property, because the uniform variables used for the jumps of our random walks are i.i.d., so two sets of uniform variables supported by disjoint boxes are independent. Furthermore, in practice, we will always use mixing to control the correlation of two events; when doing so, we say that we are decoupling them.

These remarks are expressed in the following statement, whose proof is straightforward (see Figure \ref{D:f:mixing} for an illustration). We say that an event $A$ is measurable with respect to a box $B$ if $\mathbf{1}_A$ is a measurable function of $\omega\vert_{B\cap\LLL}$ and $\{U_n^x,\,(x,n)\in B\cap\LLL\}$.

\begin{fact}[decoupling]\label{D:p:mixing}
Suppose that $\mathcal{D}_{\NN}(\ucD{D:c:mixing},2R+3,1,\alpha)$ is satisfied. Let $H>0.$ Let $B$ and $B'$ be two boxes satisfying
  $$\max(\mathrm{diam}(B),\mathrm{diam}(B'))\leqslant (2R+3)H,\quad \max(h(B), h(B'))\leqslant H,\quad \mathrm{sep}(B,B')\geqslant H.$$
Let $A$ resp. $A'$ be two events measurable with respect to $B$ resp. $B'$. We have \begin{align*}
\PP(A\cap A')\leqslant \PP(A)\PP(A')+ \ucD{D:c:mixing} H^{-\alpha}.
\end{align*}
\end{fact}

\begin{figure}[!h]
    \centering
    \begin{tikzpicture}[use Hobby shortcut,scale=0.4]
        \draw[fill, color=black!15!white] (-2.3,-1) rectangle (8.2,3);
        \draw[thick, right] (-2.3,-1) rectangle (8.2,3);
        \draw (5,1) node {$B'$};
        \draw[<->,thin,>=latex] (-2.5,-1) -- (-2.5,3);
        \draw[<->,thin,>=latex] (-2.3,3.2) -- (8.2,3.2);
        \draw[above] (2,3.3) node {$\leqslant (2R+3)H$};
        \draw[left] (-2.5,1) node {$\leqslant H$};
        \draw[fill, below left] (3,-1) circle (.05);
        \draw (3,-1) .. (3,1) .. (2,1.5) .. (2,3);
        \draw[fill, above] (2,3) circle (.05);
        \draw[fill, color=black!15!white] (-4,-4) rectangle (6.5,-8);
        \draw[thick, below left] (-4,-4) rectangle (6.5,-8);
        \draw (5,-6) node {$B$};
        \draw[fill, below left] (2,-8) circle (.05);
        \draw (2,-8).. (3,-6) .. (2,-5) .. (1,-4);
        \draw[fill] (1,-4) circle (.05);
        \draw[<->,thin,>=latex] (-4,-8.2) -- (6.5,-8.2);
        \draw[<->,thin,>=latex] (-4.2,-8) -- (-4.2,-4);
        \draw[below] (0.25,-8.3) node {$\leqslant (2R+3)H$};
        \draw[left] (-4.2,-6) node {$\leqslant H$};
        \draw[<->, thick,>=latex] (0,-4) -- (0,-1);
        \draw[right] (0,-2.5) node {$\geqslant H$};
    \end{tikzpicture}
    \caption{Illustration of the mixing property. Events describing respectively the two sample paths drawn inside the boxes can be decoupled using Fact \ref{D:p:mixing}.}
    \label{D:f:mixing}
\end{figure}

\subsubsection{Localization in boxes}
In order to apply Fact \ref{D:p:mixing}, we will need to localize the sample paths of the random walks in boxes. Also, we are going to need to consider simultaneously several random walks starting at the same time for Section \ref{D:s:v_-=v_+}, which is why we also define a horizontal interval $I_H(w)\subseteq\LLL$.

\begin{definition}\label{D:d:boxes} Let $H\in\RR_+^*$ and $w\in\RR\times\NN.$ We set
$$\left\{\begin{array}{l}I_H(w)=(w+[0,H)\times\{0\})\cap \LLL;\\
B_H(w)=w+[-(R+1)H,(R+2)H)\times [0,H).\end{array}\right.$$
\end{definition}

We refer to Figure \ref{D:f:boxes} for illustration. We set $I_H=I_H((0,0))$ and $B_H=B_H((0,0))$. We also define a class of functions that are possible samples for our random walks.

\begin{definition}\label{D:d:allowed_path}
Let $n_0\in\NN^*$. An allowed path on $\llbracket 0,n_0\rrbracket$ is a function $\chi:n \in \llbracket 0,n_0\rrbracket\mapsto \chi(n)=(\chi_1(n),\chi_2(n))\in\LLL$ such that
$$\forall n\in \llbracket 0,n_0-1\rrbracket,\;\left\{\begin{array}{l}|\chi_1(n+1)-\chi_1(n)|\leqslant R;\\
\chi_2(n+1)=\chi_2(n)+1.
\end{array}\right.$$
\end{definition}

We will be able to localize events about our random walks in boxes to which we can apply Fact \ref{D:p:mixing}. More precisely, we have the following deterministic property. Recall notation $\ell$ from Definition \ref{D:d:random_walk_discrete}.
\begin{fact}\label{D:f:allowed_path_box}
For every $H\geqslant \ell$, $w\in\RR\times\NN$ and $\chi=(\chi(n))_{n\in[0,H)\cap\NN}$ an allowed path starting in $I_H(w),$ we have
$$\{\chi(n),\,n\in[0,H)\cap\NN\}+[-\ell,\ell]\times\{0\}\subseteq B_H(w).$$
\end{fact}
We can see here that the choice of the size parameters for the boxes could be improved, since $\ell$ is a fixed parameter. But again, in practise, size parameters are not an obstacle to showing mixing properties for classical environments.

We will now work with $H\geqslant \ell.$ Mind that in Section \ref{D:s:limiting_speeds}, $H$ will be in $\NN$. The reason for writing Definition \ref{D:d:boxes} with a real $H$ is because we will need it in Section \ref{D:ss:traps}.

\subsubsection{Multi-scale renormalization}\label{D:sss:renormalization}
The proofs of several major propositions in the rest of the paper are based on the fundamental idea of multi-scale renormalization, which gives a practical method for using Fact \ref{D:p:mixing}. We now give a general idea of how such a proof works, and we will often refer to it in the future.

Assume that {\color{black}$\mathcal{D}_{\NN}(\ucD{D:c:mixing},2R+3,1,\alpha)$} holds. Suppose we want to show an estimate for the probability of some family of "bad" events $(A_H)_{H\in\NN}.$

Step 1. We start by focusing on some subsequence $(A_{H_k})_{k}$. We set $p_{k}=\PP(A_{H_k})$. We show that when $A_{H_{k+1}}$ occurs, we can find two boxes $B,B'$ satisfying
  $$\max(\mathrm{diam}(B),\mathrm{diam}(B'))\leqslant (2R+3)R_k,\quad \max(h(B), h(B'))\leqslant R_k,\quad \mathrm{sep}(B,B')\geqslant R_k$$
for some $R_k>0$, such that two events of probability $p_k$ supported by $B$ and $B'$ occur. There can be entropy for where theses two boxes are; let us say that there are $C_k$ possible couples of boxes.

Step 2. We deduce the desired estimate for $(p_k)_{k\in\N}$:
\begin{itemize}
\item Using Fact \ref{D:p:mixing} and a union bound, we get $p_{k+1}\leqslant C_k\,(p_k^2+\ucD{D:c:mixing} R_k^{-\alpha}).$
\item From this inequality we deduce the desired estimate of $p_k$ by induction on $k$. For this to work, the scales $(H_k)_k$ and the estimate to show have to be chosen properly, and some control on $C_k$ is necessary.
\end{itemize}

Step 3. We conclude using an interpolation argument to get a result for every $H\in\NN$ instead of only a subsequence.

In order to accommodate to polynomial mixing, we will use the following sequence, just as in \cite{BHT}:
\begin{align}\label{D:e:scalesL_k}\left\{\begin{array}{l}
L_0=\max(10^{10},\ell)\\
\forall k\geqslant 0,\;L_{k+1}=l_k\,L_k,\;\text{where}\;l_k=\lfloor L_k^{1/4}\rfloor.\end{array}\right.
\end{align}
{\color{black}The sequences $(H_k)$ and $(R_k)$ mentioned above will often be built using the sequence $(L_k)$.}

\subsubsection{Road map of the proof}\label{sss:roadmap}
The rest of the paper is dedicated to proving Theorem \ref{D:t:LLNd}. The road map is the same as in \cite{BHT}.

In the first part (Section \ref{D:s:limiting_speeds}), we introduce two limiting speeds $v_-,v_+\in[-R,R]$ that somehow bound the asymptotic behavior of the random walk. The major result on $v_-$ and $v_+$, Lemma \ref{D:l:bounds_on_p}, says that the probability that the random walk has average speed below $v_-$ or beyond $v_+$ between times $0$ and $H$ goes to $0$ fast enough as $H$ goes to infinity. It is the analogue of Lemma 3.1 from \cite{BHT} and its proof, given in Section \ref{D:ss:deviation_proof}, is merely an adaptation in the discrete-time finite-range setting of the proof from \cite{BHT}.

The second part of the proof (Section \ref{D:s:v_-=v_+}), just as in \cite{BHT}, consists in showing that actually $v_-=v_+$, which will give our limiting speed $\nu=v_-=v_+$. However, here the proof differs from \cite{BHT}. Indeed, {\color{black}the proof in \cite{BHT}} relies crucially on the monotonic behavior described in Section \ref{D:ss:comparison}. We will therefore need to adapt the strategy of proof and the notion of trap from \cite{BHT}. The notion of threat introduced in Definition \ref{D:d:threats}, as well as the proof of Proposition \ref{D:p:threats_proba}, however, are straightforward adaptations of \cite{BHT}. Our new ideas are particularly used in the final step of the proof, Section \ref{D:ss:delays}.

\section{Limiting speeds}\label{D:s:limiting_speeds}

\subsection{Definitions and main results}
For $y\in \LLL$ and $H\in\NN^*$, we define $V_{H}^y=\frac{1}{H}\left(X^y_{H}-\pi_1(y)\right)$ the average speed of random walk $Z^y$ between times $0$ and $H$. For $w\in\RR\times\NN$ {\color{black} and $v\in\R$}, we define
    \begin{align}\label{D:d:Ap}\left\{\begin{array}{lcl}
    A_{H,w}(v)=\left\{\exists\,y\in I_H(w),\;V_{H}^y\geqslant v\right\}
    &\text{and}&\tilde{A}_{H,w}(v)=\left\{\exists\,y\in I_H(w),\;V_{H}^y\leqslant v\right\};\\
    A_H(v)=A_{H,(0,0)}(v)
    &\text{and}&\tilde{A}_{H}(v)=\tilde{A}_{H,(0,0)}(v);\\
    p_H(v)=\PP(A_H(v))
    &\text{and}&\tilde{p}_H(v)=\PP(\tilde{A}_H(v)).
    \end{array}\right.
    \end{align}
See Figure \ref{D:f:boxes} for an illustration of the notations. Note that $\PP(A_{H,w}(v))$ and $\PP(\tilde{A}_{H,w}(v))$ actually do not depend on $w$. Indeed, $H$ being in $\NN$ here, $I_H(w)$ always contains $H$ points, so we can apply Assumption \ref{D:a:invariance}. This is why we only compute the probability of $A_H(v)$ and $\tilde{A}_H(v)$. If $H$ was not in $\NN$, $I_H(w)$ could have two different cardinalities depending on $w$, namely $\lfloor H\rfloor$ and $\lceil H\rceil$; see equations \eqref{D:d:N_k} where this problem arises.

\begin{figure}[h]
  \centering
    \begin{tikzpicture}[use Hobby shortcut]
      \draw[above left] (0, 0) rectangle (10, 2) node {$B_H(w)$};
      \draw[<->,>=latex](-0.2, 0) -- (-0.2, 2);
      \draw[left] (-0.2,1) node {$H$};
      \draw[<->,>=latex](0, -0.6) -- (4,-0.6);
      \draw[below] (2,-0.6) node {$(R+1)H$};
      \draw[<->,>=latex](6, -0.6) -- (10,-0.6);
      \draw[below] (8,-0.6) node {$(R+1)H$};
        \draw[<->,>=latex](6, -0.6) -- (4,-0.6);
      \draw[below] (5,-0.6) node {$H$};
      \draw[very thick,fill,below left] (6, 0) -- (4,0) circle (.05) node {$w$};
      \draw[very thick, above] (4, -.15) -- (4, .15) node {$I_H(w)$};
      \draw[very thick] (6, -.15) -- (6, .15);
      \draw (4.5, 0) .. (4.8, .4) .. (5.5, .8) .. (5.5, 1.2) .. (6, 1.6) .. (6.5, 2);
      \draw[dashed] (4.5, 0) -- (6.2, 2);
      \draw[fill, below right] (4.5, 0) circle (.05) node {$y = Z^y_0$};
      \draw[fill, above right] (6.5, 2) circle (.05) node {$Z^y_{H}$};
      \draw[fill, above left] (6.2, 2) circle (.05) node {$y + H (v, 1)$};
      \draw[right] (10, 1);
    \end{tikzpicture}
  \caption{Illustration of $A_{H,w}(v)$.
  Starting from the point $y\in I_H(w)$, the random walk has an average speed at least $v$ between times $0$ and $H$.}
  \label{D:f:boxes}
\end{figure}

The limiting speeds are defined by:
\begin{align*}
    &v_+=\inf\left\{v\in\RR,\;\liminf_{H\to\infty} p_H(v)=0\right\};\\
    &v_-=\sup\left\{v\in\RR,\;\liminf_{H\to\infty} \tilde{p}_H(v)=0\right\}.
\end{align*}

\begin{remark}\label{D:r:bounds_speeds}
It will be useful to use bounds on $v_+$ and $v_-$. Note that for now it is not clear that $v_-\leqslant v_+$, although it seems quite natural. Actually it will be a straightforward consequence of Lemma \ref{D:l:bounds_on_p}. But for now, we do have
\begin{align}\label{D:e:bounds_on_v_+}
\left\{\begin{array}{l}
-R\leqslant v_+\leqslant R;\\
-R\leqslant v_-\leqslant R.
\end{array}\right.
\end{align}
For instance, the first line in \eqref{D:e:bounds_on_v_+} is the consequence of the two following observations:
$$\begin{array}{ll}
\forall v<-R,&p_H(v)=\PP(\exists\, y\in I_H,\,V_H^y\geqslant v)\geqslant \PP(X_H\geqslant -RH)=1.\\
\forall v>R,&p_H(v)=\PP(\exists\,y\in I_H,\,V_H^y\geqslant v)\leqslant \PP(\exists\,y\in I_H,\,V_H^y>R)=0.
\end{array}$$
\end{remark}

\begin{remark}
It may sound unclear why we use liminfs in the definition of $v_-$ and $v_+$ instead of limsups (especially as \eqref{D:e:bounds_on_v_+} would still hold). In fact, this will be instrumental in the second part of the proof, see Section \ref{D:ss:setup}.
\end{remark}

\begin{remark}\label{D:r:monotonicity_v_+}
Note that the function defined by $v\mapsto p_H(v)$ is non-increasing. Therefore, for $v>v_+$, we have $\displaystyle\liminf_{H\to\infty} p_H(v)=0.$ Similarly, for $v<v_-$, $\displaystyle\liminf_{H\to\infty} \tilde{p}_H(v)=0$. Still, the definitions that we gave for $v_-$ and $v_+$, using liminfs, are quite weak. Yet we will be able to prove that for $v>v_+$ resp. $v<v_-$, the liminfs are actual limits, and we will even prove a precise estimate for $p_H(v)$ resp. $\tilde{p}_H(v)$ when $H$ goes to infinity.
\end{remark}

The proof of Theorem \ref{D:t:LLNd} consists of two main parts. The first part, which is the object of this section, consists in showing deviation estimates resembling \eqref{D:LLN_proba} for speeds that are below $v_-$ or beyond $v_+$: see Lemma \ref{D:l:bounds_on_p}. The second part, which corresponds to Section \ref{D:s:v_-=v_+}, aims at proving that $v_-=v_+$, which will yield the desired speed $\nu$: see Lemma \ref{D:l:v_+=v_-}.

\ncD{D:c:deviation}
\begin{lemma}\label{D:l:bounds_on_p}
Assume that Assumptions \ref{D:a:invariance} and \ref{D:a:uniform_ellipticity} are satisfied, as well as $\mathcal{D}_{\NN}(\ucD{D:c:mixing},2R+3,1,\alpha)$ for some $\ucD{D:c:mixing}>0$ and $\alpha>11$. Then, for any $\varepsilon>0$, there exists $\ucD{D:c:deviation}=\ucD{D:c:deviation}(\varepsilon)>0$ such that for every $H\in\NN^*$,
\begin{align}\label{D:deviation bounds_estimates}
    p_H(v_++\varepsilon)\leqslant \ucD{D:c:deviation} \,H^{-\alpha/4};\\
    \tilde{p}_H(v_--\varepsilon)\leqslant \ucD{D:c:deviation} \,H^{-\alpha/4}.\label{D:deviation bounds_estimates2}
\end{align}
\end{lemma}

\begin{corollary}\label{D:c:v_-<v_+} Under the same assumptions as in Lemma \ref{D:l:bounds_on_p}, we have $v_-\leqslant v_+.$
\end{corollary}
\begin{proof}
    We now prove that this is indeed a corollary of Lemma \ref{D:l:bounds_on_p}. We argue by contradiction and assume that $v_->v_+$. Let $v=\frac{v_++v_-}{2}\in (v_+,v_-)$ and $\varepsilon=\frac{v_--v_+}{2}>0$. Using \eqref{D:deviation bounds_estimates} and \eqref{D:deviation bounds_estimates2}, we have, for every $H\in\NN^*$, $p_H(v)\leqslant \ucD{D:c:deviation}(\varepsilon) H^{-\alpha/4}$ and $\tilde{p}_H(v)\leqslant \ucD{D:c:deviation}(\varepsilon) H^{-\alpha/4}$. Therefore, $1\leqslant p_H(v)+\tilde{p}_H(v)\leqslant 2\ucD{D:c:deviation}(\varepsilon) H^{-\alpha/4}\xrightarrow[]{H\to\infty}0$, which is a contradiction.
\end{proof}

This corollary can be improved as follows.
\begin{lemma}\label{D:l:v_+=v_-}
Under the same assumptions as in Lemma \ref{D:l:bounds_on_p}, we have $v_-=v_+$.
\end{lemma}
Let us now prove our main theorem as a corollary of Lemmas \ref{D:l:bounds_on_p} and \ref{D:l:v_+=v_-}. Lemma \ref{D:l:bounds_on_p} will be proved in Section \ref{D:ss:deviation_proof} and Lemma \ref{D:l:v_+=v_-} in Section \ref{D:s:v_-=v_+}.

\subsection{Proof of Theorems \ref{D:t:LLNd} and \ref{D:t:LLNc}}\label{D:ss:final_proof}

\subsubsection{Discrete time}\label{D:final_proof}
Using Lemma \ref{D:l:v_+=v_-}, we let $\nu=v_-=v_+.$ Because of Remark \ref{D:r:bounds_speeds}, $\nu\in [-R,R]$. Let $\varepsilon>0$. Using Lemma \ref{D:l:bounds_on_p}, we have, for every $n\in\NN^*$,
$$\PP\left(\left\vert\frac{X_n}{n}-\nu\right\vert\geqslant\varepsilon\right)\leqslant p_n(\nu+\varepsilon)+\Tilde{p}_n(\nu-\varepsilon)\leqslant  2\, \ucD{D:c:deviation}(\varepsilon)\,n^{-\alpha/4}.$$
Since $\alpha>4,$ $\sum_{n\in\NN^*} \PP\left(\left\vert\frac{X_n}{n}-\nu\right\vert\geqslant\varepsilon\right)<\infty$, so by the Borel-Cantelli lemma, $\frac{X_n}{n}\xrightarrow[n\to\infty]{a.s.} \nu,$ and we do have the polynomial rate of convergence \eqref{D:LLN_proba} with $\ucD{D:c:LLN}=2\ucD{D:c:deviation}.$ This concludes the proof of Theorem \ref{D:t:LLNd} when $\mathbb{T}=\NN$.

\subsubsection{Continuous time}\label{D:sss:final_proof_continuous}
In this section, we assume Theorem \ref{D:t:LLNd} and derive from it Theorem \ref{D:t:LLNc}. We first introduce some concentration bounds for the times $T_n$ used in the continuous-time setting. Recall that $(T_n)_{n\in\NN^*}$ is a Poisson process in $\RR_+^*$ of parameter $1$. This implies that:
\begin{itemize}
    \item The number $\mathcal{N}_{t,s}$ of times $T_n$ located in $[t,s]$ or $[s,t]$ (depending on whether $t\leqslant s$ or $t>s$) is a Poisson random variable of parameter $|t-s|$.
    \item The distribution of $T_n$ is a Gamma law of parameter $(n,1)$ (that is the law of density $x>0\mapsto x^{n-1}e^{-x}/\Gamma(n)$), that satisfies the following property, obtained with a Chernoff bound.
\end{itemize}
\vspace{-5pt}
\ncD{D:c:Gamma}
\begin{lemma}\label{D:l:Gamma}
For every $\varepsilon\in (0,1)$, there exists $\ucD{D:c:Gamma}=\ucD{D:c:Gamma}(\varepsilon)>0$ such that for every $n\in\NN^*$, $$\PP(\left\vert T_n-n\right\vert \geqslant \varepsilon n)\leqslant e^{-\ucD{D:c:Gamma}n}.$$
\end{lemma}

Let us consider a continuous-time environment 
$\omega$ satisfying Assumption \ref{D:a:invariance}, {\color{black}Assumption \ref{D:a:uniform_ellipticity}} and $\mathcal{D}_{\RR_+}(\ucD{D:c:mixing},a,b,\alpha)$ with $\ucD{D:c:mixing}>0$, $a>2R+3$, $b>1$ and $\alpha>11$. Recall the definition of the sequence of {\color{black}jump} times $(T_n)_{n\in\NN^*}$ from Definition \ref{D:d:random_walk_continuous}, and set $T_0=0$. We define a discrete-time environment and random walk by setting, for $n\in\NN$,
$$\left\{\begin{array}{l}
     \tilde{\omega}_n=\omega_{T_{n+1}};\\
     \tilde{X}_n=X_{T_{n}}.
\end{array}\right.$$
Now, $\tilde{X}$ has the same law as the discrete-time random walk starting at $y=(0,0)$ from Definition \ref{D:d:random_walk_discrete} on environment $\tilde{\omega}$ with the same {\color{black}jump} function $g$. This comes from the definitions and the fact that the $(U_n^x)_{n\in\NN,x\in\ZZ}$ are independent of $\omega$ and $(T_n)_{n\in\NN^*}$, so they are independent of $\tilde{\omega}.$ We now need to check that $\tilde{\omega}$ satisfies the assumptions needed for the LLN.

\begin{lemma}[Translation invariance]
    Environment $\tilde{\omega}$ satisfies Assumption \ref{D:a:invariance}.
\end{lemma}

\begin{proof}
Let $m\in\NN^*$ and $(x_0,n_0),\ldots, (x_m,n_m)\in\mathbb{L}$. Let $f$ be a non-negative measurable function on $S^m$. Then we have
    \begin{align*}
        &\EE\left[f\left(\tilde{\omega}_{n_0+n_1}(x_0+x_1),\ldots, \tilde{\omega}_{n_0+n_m}(x_0+x_m)\right)\right]\\
        &=\EE\left[f\left(\omega_{T_{n_0+n_1+1}}(x_0+x_1),\ldots, \omega_{T_{n_0+n_m+1}}(x_0+x_m)\right)\right]\\
        &=\int \EE\left[f(\omega_{t_0+t_1}(x_0+x_1),\ldots,\omega_{t_0+t_m}(x_0+x_m))\right]\, \dd\PP_{(T_{n_0},T_{n_0+n_1+1}-T_{n_0},\ldots,T_{n_0+n_m+1}-T_{n_0})}(t_0,t_1,\ldots,t_m)\\
        &=\int \EE\left[f(\omega_{t_1}(x_1),\ldots,\omega_{t_m}(x_m))\right]\,\dd\PP_{(T_{n_0+n_1+1}-T_{n_0},\ldots,T_{n_0+n_m+1}-T_{n_0})}(t_1,\ldots,t_m)\\
        &=\int \EE\left[f(\omega_{t_1}(x_1),\ldots,\omega_{t_m}(x_m))\right]\;\dd\PP_{(T_{n_1+1},\ldots,T_{n_m+1})}(t_1,\ldots,t_m)\\
        &=\EE\left[f\left(\tilde{\omega}_{n_1}(x_1),\ldots, \tilde{\omega}_{n_m}(x_m)\right)\right].
    \end{align*}
In the first equality, we simply used the definition of $\tilde{\omega}$. In the second one, we used that $(T_n)$ and $\omega$ are independent. In the third one, we used that $\omega$ satisfies Assumption \ref{D:a:invariance}. In the second-to-last equality, we used that $(T_n)$ is a Poisson process. For the last equality, we applied the same reasoning as in the lines before.
\end{proof}

\ncD{D:c:mixing2}
\begin{lemma}[Mixing]
    There exists $\ucD{D:c:mixing2}>0$ such that $\tilde{\omega}$ satisfies $\mathcal{D}_{\NN}(\ucD{D:c:mixing2},2R+3,1,\alpha)$.
\end{lemma}

\begin{proof}
Let $H,B,B',f_1,f_2$ be as in the definition of property $\mathcal{D}_{\NN}(\ucD{D:c:mixing2},2R+3,1,\alpha)$, {\color{black} where $\ucD{D:c:mixing2}$ will be chosen later}. Let $n_0,n_1,n_2,n_3\in\NN$ be such that $\pi_2(B)\cap\NN=\llbracket n_0-1,n_1-1\rrbracket$ and $\pi_2(B')\cap\NN=\llbracket n_2-1,n_3-1\rrbracket.$ We can assume without loss of generality that $n_1<n_2$. By assumption on $B$ and $B'$ we have $n_2-n_1\geqslant H$, $n_1-n_0\leqslant H$ and $n_3-n_2\leqslant H$.

Let $\mathbf{T}=(T_{n+1})_{n\in\NN}$. Note that $f_1(\tilde{\omega})$ and $f_2(\tilde{\omega})$ can be rewritten as $g_1(\omega,\mathbf{T}_{B})$ and $g_2(\omega,\mathbf{T}_{B'})$, where $\mathbf{T}_B=(T_{n+1})_{n\in\pi_2(B)}$ and similarly for $B'$. If $\mathbf{t}=(t_{n+1})_{n\in\NN}$ is an increasing sequence of $\RR_+^*$ and $\mathbf{t}_B=(t_{n+1})_{n\in \pi_2(B)},$ $g_1(\omega,\mathbf{t}_{B})$ and $g_2(\omega,\mathbf{t}_{B'})$ are respectively measurable with respect to $(\omega_{t_{n+1}}(x))_{(x,n)\in B}$ and $(\omega_{t_{n+1}}(x))_{(x,n)\in B'}$. Now, let $\varepsilon>0$ and
$$A_\varepsilon=\left\{T_{n_2}-T_{n_1}> (1-\varepsilon)H,\, T_{n_1}-T_{n_0}< (1+\varepsilon)H,\, T_{n_3}-T_{n_2}< (1+\varepsilon)H\right\}.$$
First note that using Lemma \ref{D:l:Gamma}, $$\PP(T_{n_2}-T_{n_1}\leqslant (1-\varepsilon)H)=\PP(T_{n_2-n_1}\leqslant (1-\varepsilon)H)\leqslant \PP(T_{\lceil H\rceil}\leqslant (1-\varepsilon)\lceil H\rceil)\leqslant  e^{-\ucD{D:c:Gamma}H}.$$
Let us now give an estimate of $\PP(T_{n_1}-T_{n_0}\geqslant (1+\varepsilon)H)=\PP(T_{n_1-n_0}\geqslant (1+\varepsilon)H)$. If $n_1-n_0\geqslant H/2$, we have, using Lemma \ref{D:l:Gamma},
\begin{align*}
    \PP(T_{n_1-n_0}\geqslant (1+\varepsilon)H)
    \leqslant \PP(T_{n_1-n_0}\geqslant (1+\varepsilon)(n_1-n_0))
    \leqslant e^{-\ucD{D:c:Gamma}(n_1-n_0)}\leqslant e^{-\ucD{D:c:Gamma}H/2}.
\end{align*}
Otherwise, if $n_1-n_0<H/2$, we have $\PP(T_{n_1-n_0}\geqslant (1+\varepsilon)H)\leqslant \PP(T_{\lceil H/2\rceil}\geqslant (1+\varepsilon)\lceil H/2\rceil)\leqslant e^{-\ucD{D:c:Gamma} H/2}.$ In both cases, we get the same bound. The same holds for $\PP(T_{n_3}-T_{n_2}\geqslant (1+\varepsilon)H).$ In the end, we have
    \begin{align}\label{D:e:Glaze}
        \PP(A_\varepsilon^c)\leqslant e^{-\ucD{D:c:Gamma}H}+2 e^{-\ucD{D:c:Gamma} H/2}\leqslant 3e^{-\ucD{D:c:Gamma} H/2}.
    \end{align}
     Now, on $A_\varepsilon$, $(x,T_{n+1})_{(x,n)\in B}$ and $(x,T_{n+1})_{(x,n)\in B'}$ are located in boxes $\bar{B}$ and $\bar{B}'$ (that depend on $\mathbf{T}$) such that $\mathrm{sep}(\bar{B},\bar{B}')\geqslant (1-\varepsilon)H$ and
    $$\left\{\begin{array}{l}\max(\mathrm{diam}(\bar{B}),\mathrm{diam}(\bar{B}'))\leqslant (2R+3)H\leqslant a(1-\varepsilon)H;\\
    \max(h(\bar{B}), h(\bar{B}'))\leqslant (1+\varepsilon)H\leqslant b(1-\varepsilon)H,
    \end{array}\right.$$
    provided that $\varepsilon$ is small enough (recall that $a>2R+3$ and $b>1$). This allows us to use $\mathcal{D}_{\RR_+}(\ucD{D:c:mixing},a,b,\alpha)$ in the following way: if we rewrite $A_\varepsilon$ as $\{\mathbf{T}\in \mathcal{A}_\varepsilon\}$ (with $\PP_\mathbf{T}(\mathcal{A}_\varepsilon)=\PP(A_\varepsilon)$), then for all $\mathbf{t}\in \mathcal{A}_{\varepsilon}$, we have
    \begin{align}\label{D:e:decoup}\EE[g_1(\omega,\mathbf{t}_{B})\,g_2(\omega,\mathbf{t}_{B'})]\leqslant \EE[g_1(\omega,\mathbf{t}_{B})]\,\EE[g_2(\omega,\mathbf{t}_{B'})]+\ucD{D:c:mixing}((1-\varepsilon)H)^{-\alpha}.\end{align}
Therefore, we have
\begin{align*}
    \EE[&f_1(\tilde\omega)\,f_2(\tilde\omega)]    =\EE[g_1(\omega,\mathbf{T}_{B})\,g_2(\omega,\mathbf{T}_{B'})]\\
    &=\int \,\EE[g_1(\omega,\mathbf{t}_{B})\,g_2(\omega,\mathbf{t}_{B'})]\,\dd\PP_\mathbf{T}(\mathbf{t})&\mbox{because $\omega$ and $\mathbf{T}$ are independent}\\
    &\leqslant \int \EE[g_1(\omega,\mathbf{t}_{B})]\,\EE[g_2(\omega,\mathbf{t}_{B'})]\,\dd\PP_\mathbf{T}(\mathbf{t}) + \ucD{D:c:mixing} ((1-\varepsilon)H)^{-\alpha}+3e^{-cH} &\mbox{using \eqref{D:e:Glaze} and \eqref{D:e:decoup}}
\end{align*}
Now, if we set $\theta^{t}\omega_n(x)=\omega_{t+n}(x)$ for any $t,n\in\NN$ and $x\in\ZZ$, we have
\begin{align*}
    &\int \EE[g_1(\omega,\mathbf{t}_{B})]\,\EE[g_2(\omega,\mathbf{t}_{B'})] \,\dd\PP_\mathbf{T}(\mathbf{t})\\
    &= \int \EE[g_1(\omega,\mathbf{t}_{B})]\,\EE[g_2(\theta^{t_{n_1}}\omega,\mathbf{t}_{B'}-t_{n_1})]\,\dd\PP_\mathbf{T}(\mathbf{t})\\
    &=\int \EE[g_1(\omega,\mathbf{t}_{B})]\,\EE[g_2(\omega,\mathbf{t}_{B'}-t_{n_1})] \,\dd\PP_\mathbf{T}(\mathbf{t})&\mbox{using Assumption \ref{D:a:invariance}}\\
    &=\EE[g_1(\omega,\mathbf{T}_{B})]\,\EE[g_2(\omega,\mathbf{T}_{B'}-T_{n_1})]&\mbox{using that $\mathbf{T}_{B}$ and $\mathbf{T}_{B'}-T_{n_1}$ are independent}\\
    &=\EE[g_1(\omega,\mathbf{T}_{B})]\,\EE[g_2(\omega,\mathbf{T}_{B'})]&\mbox{using the previous equalities with $g_1=1$}\\    &=\EE[f_1(\tilde\omega)]\,\EE[f_2(\tilde\omega)].
\end{align*}
In the end, if we fix a small enough $\varepsilon$, we get, for a wise choice of $\ucD{D:c:mixing2}>0,$ $\mathrm{Cov}(f_1(\tilde{\omega}),f_2(\tilde{\omega}))\leqslant \ucD{D:c:mixing2} H^{-\alpha}.$
\end{proof}

\begin{proof}[Proof of Theorem \ref{D:t:LLNc}.]
We showed that the discrete-time random walk $\tilde{X}$ satisfies the assumptions of Theorem \ref{D:t:LLNd}. As a result, there exists $\nu\in [-R,R]$ such that for every $\varepsilon>0$, there exists $\ucD{D:c:LLN}>0$ such that \begin{align}\label{D:concentration}\forall n\in\NN^*, \;\;\PP\left(\left\vert \frac{X_{T_n}}{n}-\nu\right\vert\geqslant \varepsilon\right)\leqslant \ucD{D:c:LLN}\,n^{-\alpha/4}.\end{align}
Now, let $\varepsilon>0$ and $t\geqslant 1$.
Let $n=\lfloor t\rfloor \in\NN^*.$ We have
\begin{align}\label{D:e:3terms}
\left\vert \frac{X_t}{t}-\nu\right\vert
&\leqslant \left\vert \frac{X_{T_n}}{n}-\nu\right\vert +\frac{1}{n}\left\vert X_n-X_{T_n}\right\vert+\left\vert \frac{X_n}{n}-\frac{X_t}{t}\right\vert.
\end{align}
The first term in the right-hand side of \eqref{D:e:3terms} can be bounded using \eqref{D:concentration}. As for the second term, using a fixed $\varepsilon\in (0,1)$, we have
{\color{black}
$$
\mathbb{P}\left(\frac{1}{n}\left|X_n-X_{T_n}\right| \geqslant \frac{\varepsilon}{3}\right) \leqslant \mathbb{P}\left(\mathcal{N}_{n, T_n} \geqslant \frac{n \varepsilon}{3 R}\right) .
$$
Now, fix $\xi \in(0,1)$.
\begin{align}
\nonumber\mathbb{P}\left(\mathcal{N}_{n, T_n} \geqslant \frac{n \varepsilon}{3 R}\right) & \leqslant \mathbb{P}\left(\mathcal{N}_{n, T_n} \geqslant \frac{n \varepsilon}{3 R}, T_n \in[(1-\xi) n,(1+\xi) n]\right)+\mathbb{P}\left(T_n \notin[(1-\xi) n,(1+\xi) n]\right) \\
& \leqslant \mathbb{P}\left(\mathcal{N}_{(1-\xi) n,(1+\xi) n} \geqslant \frac{n \varepsilon}{3 R}\right)+\mathbb{P}\left(T_n \notin[(1-\xi) n,(1+\xi) n]\right) \nonumber\\
& \leqslant c n^{-\alpha / 4},\label{D:concentration2}
\end{align}
using a Chernoff bound for the first term (provided that $\xi$ is small enough) and Lemma \ref{D:l:Gamma} with $T_n$ for the second one.} As for the third term in the right-hand side of \eqref{D:e:3terms}, first note that since $X$ only jumps at times $(T_n)$, we have
\begin{align*}
    \left\vert \frac{X_n}{n}-\frac{X_t}{t}\right\vert\leqslant \vert X_n\vert \left(\frac{1}{n}-\frac{1}{t}\right)+\frac{1}{t}\vert X_t-X_n\vert\leqslant \frac{R\mathcal{N}_{0,n}}{n^2} +\frac{R \mathcal{N}_{n,n+1}}{n}.
\end{align*}
Therefore, we have 
\begin{align}\label{D:concentration3}
\PP\left(\left\vert\frac{X_n}{n}-\frac{X_t}{t}\right\vert\geqslant \frac{\varepsilon}{3}\right)\leqslant \PP\left(\frac{R\mathcal{N}_{0,n}}{n^2}\geqslant \frac{\varepsilon}{6}\right)+ \PP\left(\frac{R \mathcal{N}_{n,n+1}}{n}\geqslant \frac{\varepsilon}{6}\right)\leqslant 2e^{-cn}\leqslant cn^{-\alpha/4},
\end{align}
using two Chernoff bounds in the second-to-last inequality. Putting together the probability estimates for the three terms in the right-hand side of \eqref{D:e:3terms}, we get the estimate {\color{black}in Theorem \ref{D:t:LLNc}} (with a constant that is different from the discrete case).

As for the proof of the almost sure convergence stated in Theorem \ref{D:t:LLNd}, we use \eqref{D:e:3terms} and show that all three terms in the right-hand side converge almost surely to $0$, for example, using the Borel-Cantelli lemma with \eqref{D:concentration}, \eqref{D:concentration2} and \eqref{D:concentration3} together with the fact that $\alpha>4.$
\end{proof}

\subsection{Deviations bounds: proof of Lemma \ref{D:l:bounds_on_p}}\label{D:ss:deviation_proof}
This proof is actually a direct adaptation of the continuous-time $R=1$ case that was studied in \cite{BHT}. For now, there is no extra difficulty in working in a discrete-time setting with range $R\geqslant 1$.

We will only show \eqref{D:deviation bounds_estimates}; \eqref{D:deviation bounds_estimates2} will follow using an argument of symmetry, which we state here as a remark that we will use again later on. 

\begin{remark}\label{D:r:symmetry}
Let $y=(x_0,n_0)\in \mathbb{L}$. We set $\Bar{y}=(-x_0,n_0)$, and we define the symmetrical random walk $\bar{Z}^{\bar{y}}$ using the same formulas as those defining $Z^y$, after replacing $y$ by $\bar{y}$, $g$ by $$\bar{g}:(\sigma_{-\ell},\ldots,\sigma_{\ell},u)\mapsto -g(\sigma_{\ell},\ldots,\sigma_{-\ell},u),$$ $\omega$ by $\bar{\omega}:(x,n)\mapsto \omega_n(-x)$ and $(U^x_n)_{(n,x)}$ by $(\Bar{U}_n^x=U^{-x}_n)_{(n,x)}$. That way, the symmetrical random walks satisfy
\begin{align}\label{D:e:formulas_symm}
\left\{\begin{array}{l}
\forall y\in\mathbb{L},\,\forall n\in\NN,\;\Bar{X}^{\bar{y}}_n=-X^y_n\\
\forall y\in\mathbb{L},\,\forall H\in\NN^*,\;\bar{V}^{\bar{y}}_H=-V^y_H;\\
\bar{v}_+=-v_-;\\
\bar{v}_-=-v_+,
\end{array}\right.\end{align}
where $\bar{V}^{\bar{y}}$ denotes the average speed of $\bar{Z}^{\bar{y}}$, and $\Bar{v}_-$ and $\Bar{v}_+$ denote the lower and upper limiting speeds for the symmetrical random walks. Note that in order to get the last two equalities, we use Assumption \ref{D:a:invariance}, which implies that for any $y\in\ZZ^2$, $Z$ and $Z^y$ have the same law. Now, $\bar{g}$, $\bar{\omega}$ and $(\Bar{U}_n^x)$ satisfy the same assumptions as $g$, $\omega$ and $(U_n^x)$, so everything we show for $X^y$ is also satisfied by $\bar{X}^{\bar{y}}$. If $J$ is a horizontal interval of $\RR\times\NN$, we denote by $\Bar{J}$ the symmetrical interval given by $\{(-x,t),(x,t)\in J\}.$
\end{remark}

Recall that $I_H=I_H((0,0))$. Using Remark \ref{D:r:symmetry}, we have
\begin{align*}
    \Tilde{p}_H(v_--\varepsilon)
    =\PP(\exists\,y\in I_H,\,V_H^y\leqslant v_--\varepsilon)
    =\PP(\exists\,y\in I_H,\,\bar{V}_H^{\bar{y}}\geqslant \Bar{v}_++\varepsilon)
    =\PP(\exists\,y\in I_H,\,\bar{V}_H^{y}\geqslant \Bar{v}_++\varepsilon),
\end{align*}
using Assumption \ref{D:a:invariance} for the last equality (note that $I_H$ and $\bar{I}_H=-I_H$ have the same number of points). Now, an upper bound for this last term is given by \eqref{D:deviation bounds_estimates}, using Remark \ref{D:r:symmetry}, which proves \eqref{D:deviation bounds_estimates2}. This allows us to show only \eqref{D:deviation bounds_estimates}.
\ncDk{D:k:first_scale}
\vspace{5pt}

Now, in order to show \eqref{D:deviation bounds_estimates}, we will follow the method given in Section \ref{D:sss:renormalization} step by step.

\subsubsection{Estimate along a subsequence}

We are first going to show an estimate similar to \eqref{D:deviation bounds_estimates} along a subsequence. More precisely, we are going to show the following result.
\begin{lemma}\label{D:l:step2}
    Let $v>v_+$. There exists $\ucDk{D:k:first_scale}=\ucDk{D:k:first_scale}(v)\in\NN$ and $h=h(v)\in [1,\infty)$ such that for every $k\geqslant \ucDk{D:k:first_scale}$, $hL_k\in\NN$ and \begin{align}\label{D:e:along_subsequence_estimate}
    p_{h L_k}(v)\leqslant L_k^{-\alpha/2}.
\end{align}
\end{lemma}

\begin{proof} We fix $v>v_+$ for the whole proof. For now, we also fix a parameter $h\in\NN^*$ that we will choose later on, which will act as a zooming parameter designed for the base case of the renormalization to work. \vspace{-10pt}
\paragraph{Step 1.} We start by linking what happens at scales $hL_k$ and $hL_{k+1}$. The key idea of the proof is that if we choose a sequence of speeds $(v_k)$ properly, on the event that there exists a random walk going faster than $v_{k+1}$ at scale $hL_{k+1}$, there exists two well-separated boxes in scale $hL_k$ on which there exist random walks going faster than $v_k$. Let $\ucDk{D:k:first_scale}=\ucDk{D:k:first_scale}(v)\in\NN$ be a constant satisfying the following technical conditions:
\begin{align}
\label{D:e:condition1}
&\sum_{k\geqslant \ucDk{D:k:first_scale}} \frac{4R}{l_k}<\frac{v-v_+}{2};\\
\label{D:e:condition2} &\forall k\geqslant \ucDk{D:k:first_scale},\;\;(2R+3)^2(\ucD{D:c:mixing}+1)\,L_k^{1-3\alpha/8}\leqslant 1.
\end{align}
Note that \eqref{D:e:condition2}, which will be instrumental in Step 2 of this proof, makes sense because $\alpha\geqslant 3$.
We define a sequence of speeds (see Figure \ref{D:fig:v_k} for illustration) by setting
\begin{align}\label{D:e:speed_sequence}\left\{
\begin{array}{l}
v_{\ucDk{D:k:first_scale}}=\displaystyle\frac{v+v_+}{2};\\
\forall k\geqslant \ucDk{D:k:first_scale},\;\;v_{k+1}=v_k+\displaystyle\frac{4R}{l_k}.
\end{array}\right.\end{align}
Using \eqref{D:e:condition1}, we have $$v_k\xrightarrow[k\to\infty]{\nearrow} v_\infty=v_{\ucDk{D:k:first_scale}}+\sum_{k\geqslant \ucDk{D:k:first_scale}} \frac{4R}{l_k}<v.$$

\begin{figure}[h]
  \begin{center}
    \begin{tikzpicture}
      \draw (-1, 0) -- (10, 0);
      \draw (0, .2) -- (0, -.2) node [below] {$v_+$};
      \draw (3, .2) -- (3, -.2) node [below] {$v_{\ucDk{D:k:first_scale}}$};
      \draw (6, .2) -- (6, -.2) node [below] {$v$};
      \draw (4, .2) -- (4, -.2) node [below] {$v_{\ucDk{D:k:first_scale}+1}$};
      \draw (4.5, .2) -- (4.5, -.2);
      \draw (4.9, .2) -- (4.9, -.2);
      \draw (4.75, .2) -- (4.75, -.2) node [below] {$\dots$};
      \draw (5.5, .2) -- (5.5, -.2) node [below] {$v_{\infty}$};
    \end{tikzpicture}
  \caption{The sequence of speeds $(v_k)_{k\geqslant \ucDk{D:k:first_scale}}$.}
  \label{D:fig:v_k}
  \end{center}
\end{figure}

Let $k\geqslant \ucDk{D:k:first_scale}$. We set \begin{align}\label{D:CCC}
    \CCC=\CCC_{h,k}=\{(ihL_k,jhL_k)\in hL_k \ZZ^2,\,-(R+1)l_k\leqslant i<(R+2)l_k,\,0\leqslant j<l_k\}.
\end{align}
This set has cardinality $(2R+3)\,l_k^2$ and satisfies \begin{align}
\bigcup_{w\in \CCC} I_{hL_k}(w)=B_{hL_{k+1}}\cap \left(\R\times hL_k\Z\right),\label{D:def_C}
\end{align}
where the union above is disjoint (recall that $B_{hL_{k+1}}=B_{hL_{k+1}}((0,0))$).

Recall definitions \eqref{D:d:Ap}. We claim that on $A_{hL_{k+1}}(v_{k+1})$, there exist two points $w_1,\,w_2\in \CCC$ such that \begin{align}\label{D:e:Glazou}
\left\{\begin{array}{l}A_{hL_k,w_1}(v_k)\text{ and }A_{hL_k,w_2}(v_k)\text{ occur;}\\ \vert\pi_2(w_1)-\pi_2(w_2)\vert\geqslant 2hL_k.\end{array}\right.
\end{align}
Indeed, let us assume that $A_{hL_{k+1}}(v_{k+1})$ occurs. It is enough to show that there exist $w_1',\,w_2',\,w_3'\in \CCC$ such that $A_{hL_k,w_i'}(v_k)$ occurs for $i=1,2,3$ and $(\pi_2(w_i'))_{i=1,2,3}$ are distinct. Assume by contradiction that this does not hold. It means that we can find $j_1,\,j_2\in\{0,\ldots,l_k-1\}$ such that 
\begin{align}\label{D:e:cascade}
\text{for all $j\notin \{j_1,\,j_2\}$, } (A_{hL_k,w}(v_k))^c\text{ occurs for all $w\in\CCC$ such that $\pi_2(w)=jhL_k$}
\end{align}
Fix $y\in I_{hL_{k+1}}$. Remark that because of Fact \ref{D:f:allowed_path_box}, for every $j\in\llbracket 0,l_k-1\rrbracket $, $Z^y_{jhL_k}\in B_{hL_{k+1}}$. Now its second coordinate is in $hL_k\ZZ$ by construction, so using \eqref{D:def_C}, $Z^y_{jhL_k}$ is in $I_{hL_k}(w)$ for some $w\in \CCC$. Therefore, by \eqref{D:e:cascade}, for $j\notin \{j_1,\,j_2\}$, we can upper-bound the displacement between times $jhL_k$ and $(j+1)hL_k$ by $v_khL_k$. For $j_1$ and $j_2$, we can simply upper-bound the displacement by $RhL_k$, by definition of the {\color{black}jump} range $R$. In the end, we get
\begin{align*}
    X_{hL_{k+1}}^y-\pi_1(y)&=\sum_{j=0}^{l_k-1} \left(X_{(j+1)hL_k}^{y}-X^y_{jhL_k}\right)\\
    &\leqslant (l_k-2)v_khL_k+2RhL_k\\
    &=v_khL_{k+1}+\frac{2(R-v_k)}{l_k}\,hL_{k+1}\\
    &<\left(v_k+\frac{4R}{l_k}\right)\,hL_{k+1}&\mbox{using that $v_k>v_+$ and \eqref{D:e:bounds_on_v_+}}\\
    &=v_{k+1}hL_{k+1}&\mbox{by definition of $v_{k+1}$},
\end{align*}
which contradicts the fact that $A_{hL_{k+1}}(v_{k+1})$ occurs. This proves \eqref{D:e:Glazou}. Now, note that in \eqref{D:e:Glazou}, events $A_{hL_k,w_1}(v_k)$ and $A_{hL_k,w_2}(v_k)$ are measurable with respect to boxes $B_{hL_k}(w_1)$ and $B_{hL_k}(w_2)$, using Fact \ref{D:f:allowed_path_box}, and these two boxes satisfy the assumptions of Fact \ref{D:p:mixing} with $H=hL_k$, by Definition \ref{D:d:boxes} and the fact that $|\pi_2(w_1)-\pi_2(w_2)|\geqslant 2hL_k$. Therefore we will be able to decouple them.\vspace{-7pt}
\paragraph{Step 2.}
We now define $h=h(v)$ and use induction on $k\geqslant \ucDk{D:k:first_scale}$ to conclude. Using Remark \ref{D:r:monotonicity_v_+}, note that since $v_{\ucDk{D:k:first_scale}}>v_+$, $$\liminf_{H\to\infty,\,H\in\NN} p_H(v_{k_0})=0.$$ Therefore there exists an integer $H\geqslant L_{\ucDk{D:k:first_scale}}$ such that $p_H(v_{\ucDk{D:k:first_scale}})\leqslant L_{\ucDk{D:k:first_scale}}^{-\alpha/2}$. Let $h=H/L_{\ucDk{D:k:first_scale}}\in [1,\infty).$ Using that $L_k$ is a multiple of $L_{\ucDk{D:k:first_scale}}$ for every $k\geqslant \ucDk{D:k:first_scale}$, we have $hL_k\in\NN$ for all $k\geqslant \ucDk{D:k:first_scale}$, and
\begin{align}\label{D:e:base_case}
p_{h L_{\ucDk{D:k:first_scale}}}(v_{k_0})\leqslant L_{\ucDk{D:k:first_scale}}^{-\alpha/2}.
\end{align}
Now, note that if we show that for every $k\geqslant \ucDk{D:k:first_scale}$, we have \begin{align}\label{D:e:estimate_to_show} p_{hL_k}(v_k)\leqslant L_k^{-\alpha/2},\end{align} then we get \eqref{D:e:along_subsequence_estimate} using the monotonicity property in Remark \ref{D:r:monotonicity_v_+}. So what we will do now is show estimate \eqref{D:e:estimate_to_show} by induction on $k\geqslant \ucDk{D:k:first_scale}$. The base case is just \eqref{D:e:base_case}. Now, suppose that \eqref{D:e:estimate_to_show} is satisfied for a fixed $k\geqslant \ucDk{D:k:first_scale}$. Using our claim shown in Step 1, we have, setting $\CCC=\CCC_{h,k}$ and using a union bound, 
    \begin{align*}
        \PP\left(A_{hL_{k+1}}(v_{k+1})\right)
        &\leqslant |\CCC|^2\,\sup_{w_1,w_2\in \CCC\atop |\pi_2(w_1)-\pi_2(w_2)|\geqslant 2hL_k} \PP\left(A_{hL_k,w_1}(v_k)\cap A_{hL_k,w_2}(v_k)\right)\\
        &\leqslant (2R+3)^2\,l_k^4\,\left(p_{hL_k}(v_k)^2+\ucD{D:c:mixing}(hL_k)^{-\alpha}\right),
    \end{align*}
    using Fact \ref{D:p:mixing}. Then, using the induction assumption,
    \begin{align*}
        \frac{p_{h L_{k+1}}(v_{k+1})}{L_{k+1}^{-\alpha/2}}
        &\leqslant (2R+3)^2\,L_{k+1}^{\alpha/2}\,l_k^4\,\left(p_{hL_k}(v_k)^2+\ucD{D:c:mixing}(hL_k)^{-\alpha}\right)\\
        &\leqslant (2R+3)^2\,L_k^{5\alpha/8+1}\,(L_k^{-\alpha}+\ucD{D:c:mixing} L_k^{-\alpha})&\mbox{because $h\geqslant 1$}\\
        &\leqslant (2R+3)^2(\ucD{D:c:mixing}+1)\,L_k^{1-3\alpha/8}\\
        &\leqslant 1 &\mbox{by \eqref{D:e:condition2}}.
    \end{align*}
    This concludes the induction and thus the proof of \eqref{D:e:along_subsequence_estimate} for every $k\geqslant \ucDk{D:k:first_scale}$.
\end{proof}

\ncDk{D:k:interpolation}\ncDk{D:k:int}
\subsubsection{Interpolation} To end the proof of \eqref{D:deviation bounds_estimates}, it remains to interpolate \eqref{D:e:along_subsequence_estimate} to get something valid for any $H\in\NN^*$. Recall that constants $\ucDk{D:k:first_scale}$ and $h$ from Lemma \ref{D:l:step2} actually depend on the speed in the estimate. Let us fix $v>v_+$, and consider $v'=\frac{v_++v}{2}$, $h=h(v')$ and $\ucDk{D:k:interpolation}\geqslant \ucDk{D:k:first_scale}(v')$ such that \begin{align}\label{D:e:condition_k}
L_{\ucDk{D:k:interpolation}}^{1/10}>\frac{R}{v-v'}.
\end{align}
Let $H\geqslant(hL_{\ucDk{D:k:interpolation}})^{11/10}$ and let $\ucDk{D:k:int}\geqslant \ucDk{D:k:interpolation}$ be such that
\begin{align}\label{D:11/10}(hL_{\ucDk{D:k:int}})^{11/10}\leqslant H<(hL_{\ucDk{D:k:int}+1})^{11/10}.
\end{align}
Since $\ucDk{D:k:int}\geqslant k_0(v')$ and $v'>v_+$, \eqref{D:e:along_subsequence_estimate} ensures that
\begin{align}\label{D:cocorico}
p_{hL_{\ucDk{D:k:int}}}(v')\leqslant L_{\ucDk{D:k:int}}^{-\alpha/2}.
\end{align}
Now, consider box $B_H$. In order to link what happens at scale $H$ with information at scale $h L_{\ucDk{D:k:int}}$, similarly to \eqref{D:CCC}, we define a set
$$\bar{\CCC}=\left\{(ihL_{\ucDk{D:k:int}},jhL_{\ucDk{D:k:int}})\in hL_{\ucDk{D:k:int}} \ZZ^2,\,-(R+1)\left\lceil \frac{H}{hL_{\ucDk{D:k:int}}}\right\rceil\leqslant i<(R+2)\left\lceil \frac{H}{hL_{\ucDk{D:k:int}}}\right\rceil,\, 0\leqslant j<\lfloor H/(hL_{\ucDk{D:k:int}})\rfloor\right\},$$
which satisfies 
\begin{align}\label{D:d:barC}
\bigcup_{w\in\bar{\CCC}} I_{hL_{\ucDk{D:k:int}}}(w) \supseteq B_H\cap (\RR\times hL_{\ucDk{D:k:int}} \ZZ),
\end{align}
and whose cardinality can be bounded from above as follows:
\begin{align}\label{D:cardinality}|\bar{\CCC}|\leqslant c\left(\frac{H}{hL_{\ucDk{D:k:int}}}\right)^2\leqslant c\left(\frac{(hL_{\ucDk{D:k:int}+1})^{11/10}}{hL_{\ucDk{D:k:int}}}\right)^2\leqslant c(h)\,L_{\ucDk{D:k:int}}^{3/4}.
\end{align}
Now, let us assume that $(A_{hL_{\ucDk{D:k:int}},w}(v'))^c$ occurs for every $w\in \bar{\CCC}$. Then, using \eqref{D:d:barC}, we have, for any $y\in I_H$,
\begin{align}
    X^y_{\lfloor H/(hL_{\ucDk{D:k:int}})\rfloor\,hL_{\ucDk{D:k:int}}} -\pi_1(y)&=\sum_{j=0}^{\lfloor H/(hL_{\ucDk{D:k:int}})\rfloor-1}\left(X_{hL_{\ucDk{D:k:int}}}^{Z^y_{jhL_{\ucDk{D:k:int}}}}-X^y_{jhL_{\ucDk{D:k:int}}}\right)\tag*{}\\
    &\leqslant v'\,\lfloor H/(hL_{\ucDk{D:k:int}})\rfloor\,hL_{\ucDk{D:k:int}}\tag*{}\\
    &\leqslant v'H.\label{D:<v'H}
\end{align}
Now, between time $\lfloor H/(hL_{\ucDk{D:k:int}})\rfloor\,hL_{\ucDk{D:k:int}}$ and $H$, we have
\begin{align}\label{D:<(v-v')H}
    X^y_H-X^y_{\lfloor H/(hL_{\ucDk{D:k:int}})\rfloor\,hL_{\ucDk{D:k:int}}}\leqslant hL_{\ucDk{D:k:int}}R< (v-v')H,
\end{align}
using \eqref{D:e:condition_k}, \eqref{D:11/10} and $h\geqslant 1$. Combining \eqref{D:<v'H} and \eqref{D:<(v-v')H}, we get $X_H^y<\pi_1(y)+vH$. This being true for every $y\in I_H$, this means that $A_H(v)$ cannot occur. Therefore we have
\begin{align*}
    \PP(A_{H}(v))
    \leqslant \PP\left(\underset{w\in \bar{\CCC}}{\bigcup} A_{hL_{\ucDk{D:k:int}},w}(v')\right)
    &\leqslant |\bar{\CCC}|\,L_{\ucDk{D:k:int}}^{-\alpha/2}&\mbox{using \eqref{D:cocorico}}\\
    &\leqslant c(h)\,L_{\ucDk{D:k:int}}^{3/4}L_{\ucDk{D:k:int}}^{-\alpha/2}&\mbox{using \eqref{D:cardinality}}\\
    &\leqslant c(h)\,L_{\ucDk{D:k:int}}^{-7\alpha/20}&\mbox{using that $\alpha\geqslant 5$}\\
    &\leqslant c(h)\,H^{-\alpha/4}&\mbox{using \eqref{D:11/10}.}
\end{align*}
To conclude, it suffices to adapt $\ucD{D:c:deviation}$ in order for this estimate to hold also for small values of $H$, which proves \eqref{D:deviation bounds_estimates}.

\section{Equality of the limiting speeds: proof of Lemma \ref{D:l:v_+=v_-}}\label{D:s:v_-=v_+}
Now that we have the inequalities given by Lemma \ref{D:l:bounds_on_p}, we want to prove that $v_-=v_+$, namely Lemma \ref{D:l:v_+=v_-}. Mind that contrary to the proof of Lemma \ref{D:l:bounds_on_p}, here the proof from \cite{BHT} for the $R=1$ case cannot be adapted straight away. Indeed, the fact that random walks can jump at larger range allow them to cross paths (see Remark \ref{D:r:monotonicity_BHT} and Section \ref{D:ss:comparison}).

First recall that $v_-\leqslant v_+$, on account of Corollary \ref{D:c:v_-<v_+}. So, for the rest of this section, we argue by contradiction and assume that $v_-<v_+$. Now, the idea is to show that on some scales, the random walk's average speed will be close to (although larger than) $v_-$, so it will accumulate a delay with respect to $v_+$. Now it is hard to catch up this delay (and therefore attain an average speed close to $v_+$ later), because the probability that the random walk has average speed at least $v_++\varepsilon$ (where $\varepsilon>0$) at scale $H$ is vanishing with $H$, due to Lemma \ref{D:l:bounds_on_p}. Nonetheless, it might be able to compensate the delay by going faster than $v_++ \varepsilon(H)$ for some well-chosen $\varepsilon(H)>0$. This is why we need precise estimates on the delay, and these will be given using a notion of trap adapted from that in \cite{BHT}.

In the rest of this section, we set  
\begin{align}\label{D:d:delta}\begin{array}{ccc}
v_0=\displaystyle\frac{v_-+v_+}{2}&\text{and}&\delta=\displaystyle\frac{v_+-v_-}{6R}.\end{array}
\end{align}
Note that because of Remark \ref{D:r:bounds_speeds}, we have $0<\delta\leqslant \frac{1}{3}.$

\subsection{Setup of the proof}\label{D:ss:setup}
We will be working with three different scales $h_k\ll H_k\ll \mathbf{H}_k$ where $k\in\mathcal{S}$, an unbounded subset of $\NN$. Our goal will be to prove that there exists $\varepsilon>0$ such that
\begin{align}\label{D:e:contradiction}
p_{\mathbf{H}_k}(v_+-\varepsilon)\xrightarrow[k\to\infty\atop k\in \mathcal{S}]{}0,
\end{align}
which is a contradiction with the definition of $v_+$ (recall definitions \eqref{D:d:Ap}). Here it is crucial that in the definition of $v_+$ we used a liminf, because \eqref{D:e:contradiction} only provides a limit on the subsequence given by $(\mathbf{H}_k)_{k\in\mathcal{S}}.$ In practice, what we try to do is to delay random walks on these scales, meaning that we give an upper bound for their average speeds that is away from $v_+$. In order to do that, we adapt notions of traps and threats introduced in \cite{BHT}. However, in \cite{BHT}, traps and threats automatically imply delays on account of the fundamental monotonicity property that we discussed in Remark \ref{D:r:monotonicity_BHT}. Here, we will have to work more to get delays, using the much weaker coalescing property \eqref{D:e:coupling_property} along with Assumption \ref{D:a:uniform_ellipticity}.
\begin{itemize}[leftmargin=*]
\item We are first going to define a notion of traps: a trapped point will correspond to an area in space-time where a random walk runs an increased risk of being delayed at scale $h_k$: see Section \ref{D:ss:traps}.
\item However, the probability for a point to be trapped, although it is at least $1/2$, will not be arbitrarily close to $1$. In order to increase that probability, we will look for trapped points among a large number $r$ of points along a segment rooted at a fixed point, and we will say that this latter point is threatened when we can find at least one trapped point. We will be able to show that a point is threatened with probability close to $1$ by choosing $r$ large enough: see Proposition \ref{D:p:threats_proba}. The interest for our purpose of delaying random walks is that a random walk near a threatened point runs an increased risk of being delayed at scale $H_k=rh_k$: see Section \ref{D:ss:delays}.
\item The problem now is that although a fixed point has a high probability of being threatened, we have to make sure that a lot of the points actually encountered by our random walk are near threatened points. We will show that for a good choice of our third set of scales $(\mathbf{H}_k)$, with a high probability, half of the points encountered are indeed near a threatened point: see Section \ref{D:ss:threatened_paths}.
\item Finally, using Assumption \ref{D:a:uniform_ellipticity}, we will argue that with probability at least $\gamma$, a random walk running the risk of being delayed is actually delayed. Combined with the previous points, we will be able to get a delay at scale $\mathbf{H}_k$: see Section \ref{D:ss:delays}.
\end{itemize}
Although it is a bit premature, we refer the reader to Figure \ref{D:f:density} for an illustration of this line of reasoning.

\subsection{Trapped points}\label{D:ss:traps}
\ncDk{D:k:paths} \ncD{D:c:threats}
We define the following constants (that only depend on $\ucD{D:c:mixing}$ and $\alpha$), which will appear naturally in the proof of Proposition \ref{D:p:threats_proba}:
\begin{align}\label{D:constants_threats}
    \left\{\begin{array}{l}
        j_0=\max\left(0,\left\lceil \log_3\left(\frac{\sqrt{2}\, \ucD{D:c:mixing}}{\sqrt{2}-1}\right)\right\rceil\right);\\
        j_1=\max\left(j_0+\lceil 2(\alpha \log_2(3)+1)\rceil,\,\frac{\log_3(4\ucD{D:c:mixing})}{\alpha}+1\right);\\
        \ucD{D:c:threats}=3^{\alpha(j_1+1)}.
    \end{array}\right.
\end{align}
Recall \eqref{D:d:delta}. Let $\ucDk{D:k:paths}\in\NN$ be such that
\begin{align}
\label{D:condition_on_k0} &\delta L_{\ucDk{D:k:paths}}\geqslant 12;\\
\label{D:e:condition_on_k1} &\forall k>\ucDk{D:k:paths},\;\;(2R+3)^2\,\left(\ucD{D:c:threats}^2\,\frac{4^{2\alpha+4}}{\delta^2}+\ucD{D:c:mixing}\right) L_k^{-(3\alpha-23)/20}\leqslant \ucD{D:c:threats}\, \frac{4^{\alpha+2}}{\delta}.
\end{align}
This is possible because $\alpha\geqslant 8$. These conditions will appear in Section \ref{D:ss:threatened_paths}.

For $k\in\NN$, we set
\begin{align}\label{D:hk}
h_k=L_{\ucDk{D:k:paths}}L_k.
\end{align}
Note that for every $k\in\NN$, $\delta h_k>1$, using \eqref{D:condition_on_k0}. For $w\in\RR\times\NN$ and $k\in\NN$, we define
\begin{align}\label{D:d:N_k}
    \left\{\begin{array}{l}
    J_k(w)=(w+[\delta h_k,2\delta h_k)\times\{0\})\cap\LLL;\\
    N_k(w)=\#J_k(w);\\
    N_k^-(w)=\#\{y\in J_k(w),\,V^y_{h_k}\leqslant v_0\};\\
    N_k^+(w)=\#\{y\in J_k(w),\,V^y_{h_k}\geqslant v_0\}.
    \end{array}\right.
\end{align}
As usual, when $w$ is not specified, we mean $w=(0,0)$. Mind that depending on $w$, $\delta$ and $k$, $N_k(w)$ can be equal to $\lfloor \delta h_k\rfloor$ or $\lceil \delta h_k\rceil$, and we cannot restrict ourselves to $w\in\ZZ\times\NN$; see for instance the $w_i$ introduced in Definition \ref{D:d:threats}.

By definition, $J_k(w)$ is a horizontal interval on the right of $w$, so the idea is that the random walks starting in this interval which have small enough speeds will prevent a random walk starting near $w$ to have too large a speed, with a positive probability. This is the reason of our uniform ellipticity assumption and of our coupling (recall \eqref{D:e:coupling_property}). This heuristic idea will be made clearer in Section \ref{D:ss:delays}.

We now define a notion of traps for when a fixed proportion of random walks in $J_k(w)$ do have speed at most $v_0$.

\begin{definition}\label{D:d:traps}
We say that $w\in\RR\times\NN$ is $k$-trapped if $N_k^-(w)\geqslant N_k(w)/3$.
\end{definition}

\begin{figure}[h]
  \begin{center}
    \begin{tikzpicture}[use Hobby shortcut, scale=0.4]
        \draw[fill, below left] (0,0) circle (.05) node {$w$};
        \draw[below left] (3,0) -- (6,0);
        \draw[below] (4.5,-0.2) node {$J_k(w)$};
        \draw[fill] (3,0) circle (.05);
        \draw[fill] (3.5,0) circle (.05);
        \draw[fill] (4,0) circle (.05);
        \draw[fill] (4.5,0) circle (.05);
        \draw[fill] (5,0) circle (.05);
        \draw[fill] (5.5,0) circle (.05);
        \draw[fill] (6,0) circle (.05);
        \draw (3,-0.15) -- (3,0.15);
        \draw (6,-0.15) -- (6,0.15);
        \draw (3,0) .. (3,1) .. (2,2) .. (1.5,3) .. (1.5,4) .. (2,5) .. (1,6) .. (0,7);
        \draw (3.5,0) .. (3.5,1) .. (4,2) .. (2.5,3) .. (3,4) .. (2,5) .. (1,6) .. (1,7);
        \draw (4.5,0) .. (5.1,1) .. (3.5,2) .. (3,3) .. (4,4) .. (3,5) .. (2, 6) .. (1.5,7);
        \draw (5.5,0) .. (4.8,1) .. (5,2) .. (4,3) .. (5,4) .. (4,5) .. (3,6) .. (2,7);
        \draw[->, thick,>=latex, above] (-2,-1) -- (-2,8.3) node {time};
        \draw[left] (-1.85,0) -- (-2.15,0) node {$\pi_2(w)$};
        \draw[left] (-1.85,7) -- (-2.15,7) node {$\pi_2(w)+h_k$};
    \end{tikzpicture}
    \caption{Illustration of $w$ being $k$-trapped -- here $v_0=0$. At least a third of the random walks starting in $J_k(w)$ reach height $\pi_2(w)+h_k$ with a nonpositive average speed. When $R\geqslant 2$, sample paths can jump over each other as in the drawing.}
    \label{D:f:traps}
    \end{center}
\end{figure}

There is no reason for one third of the random walks in $J_k(w)$ to go slower than $v_0$. However, note that we have $N_k^-+N_k^+\geqslant N_k$, so if $N_k^-<N_k/2$, then $N_k^+\geqslant N_k/2$. This ensures that we have either $\PP(N_k^-\geqslant N_k/2)\geqslant 1/2$ or $\PP(N_k^+\geqslant N_k/2)\geqslant 1/2$. As a result, at least one of the two following scenarios occurs: there exists infinitely many $k\in\NN$ such that $\PP(N_k^-\geqslant N_k/2)\geqslant 1/2$; there exists infinitely many $k\in\NN$ such that $\PP(N_k^+\geqslant N_k/2)\geqslant 1/2$.

If the second scenario occurs, then the first one occurs for the symmetrical random walks defined in Remark \ref{D:r:symmetry}. Indeed, through the process of replacing $X^y$ by $\bar{X}^y$ for every $y$, using \eqref{D:e:formulas_symm}, we can see that $v_0$ becomes $\Bar{v}_0=-v_0$; as a result, if $k\in\NN$ is such that $N_k^+\geqslant N_k/2$, then $$\#\{y\in \bar{J}_k,\,\bar{V}^{y}_{h_k}\leqslant \bar{v}_0\}\geqslant \#\bar{J}_k/2.$$
Therefore, if $\PP(N_k^+\geqslant N_k/2)\geqslant 1/2$, then $\PP(\Bar{N}_k^-\geqslant N_k/2)\geqslant 1/2$, where $\Bar{N}_k^-=\#\{y\in J_k,\,\Bar{V}_{h_k}^y\leqslant \Bar{v}_0\}$ (we used translation invariance to replace $J_k$ by $\Bar{J}_k$, and they both have cardinality $N_k$). Therefore, if the second scenario occurs, then the first one occurs for the symmetrical random walks, for which $\bar{\delta}=\frac{\bar{v}_+-\bar{v}_-}{6R}=\delta$. Thus, up to showing that $\Bar{v}_+=\Bar{v}_-$ instead of $v_+=v_-$ (which is equivalent, using \eqref{D:e:formulas_symm}), we may and do assume without loss of generality that the first scenario occurs.

Now, mind that despite translation invariance, $\PP(N_k^-(w)\geqslant N_k(w)/2)$ may actually depend on $w$, for $N_k(w)$ can take two values depending on $w$.

Recall \eqref{D:condition_on_k0}. If $N_k=\lceil \delta h_k\rceil$, then for all $w\in\RR\times\NN$ such that $N_k(w)=\lfloor \delta h_k\rfloor$, we have, using translation invariance,
    \begin{align*}
    \PP(N_k^-\geqslant \lceil \delta h_k\rceil/2)
\leqslant \PP(N_k^-(w)\geqslant \lceil \delta h_k\rceil/2-1)
\leqslant \PP(N_k^-(w)\geqslant \lfloor \delta h_k\rfloor/2-1)
\leqslant \PP(N_k^-(w)\geqslant \lfloor \delta h_k\rfloor/3).
\end{align*}
If $N_k=\lfloor \delta h_k\rfloor$, then for all $w\in\RR\times\NN$ such that $N_k(w)=\lceil \delta h_k\rceil$, we have, using translation invariance,
$$\PP(N_k^-\geqslant \lfloor \delta h_k\rfloor/2)\leqslant \PP(N_k^-(w)\geqslant \lfloor \delta h_k\rfloor/2)\leqslant \PP\left(N_k^-(w)\geqslant \frac{\lceil \delta h_k\rceil -1}{2}\right)\leqslant \PP(N_k^-(w)\geqslant \lceil \delta h_k\rceil/3).$$

In both lines of inequalities, the last event is precisely $\{\text{$w$ is $k$-trapped}\}$. At the end of the day, putting together scenario 1 and these inequalities, we can conclude that there exists an unbounded subset of $\NN$, that we denote by $\mathcal{S}$, depending on $\delta$, such that
\begin{align}\label{D:e:for_traps}
\forall k\in\mathcal{S},\;\; \forall w\in\RR\times\NN,\;\;\PP(\text{$w$ is $k$-trapped})\geqslant 1/2.
\end{align}

\subsection{Threatened points}\label{D:ss:threats}
The problem with \eqref{D:e:for_traps} is that it does not guarantee that a point is trapped with high probability. This is why we now strengthen the notion of traps by introducing a notion of threats.

\begin{definition}\label{D:d:threats}
Let $k\in\NN$, $r\in\NN^*$ and $w\in\RR\times\NN$. We say that $w$ is $(k,r)$-threatened if there exists $i\in\llbracket 0,r-1\rrbracket$ such that $w_i:=w+ih_k(v_+,1)\in\RR\times\NN$ is $k$-trapped.
\end{definition}

As we said before, the interest of a threatened point is that a random walk starting in its vicinity runs a risk of being delayed at scale $H_k$. The idea is the following: since our random walk cannot go much faster than $v_+$ with high probability (using Lemma \ref{D:l:bounds_on_p}), it cannot be too much on the right of $w_i$ at time $\pi_2(w)+ih_k$. Therefore, if $w_i$ is trapped, the random walk runs the risk of being delayed between times $\pi_2(w)+ih_k$ and $\pi_2(w)+(i+1)h_k$, and this delay should be sufficient to entail a delay at scale $H_k=rh_k$. This heuristic idea will be made clearer in Section \ref{D:ss:delays}. For now, let us prove that if we choose $r$ large enough, the probability for a point to be threatened can be made arbitrarily close to $1$.

{\color{black} Recall \eqref{D:constants_threats} : $\ucD{D:c:threats}$ is a positive constant that does not depend on $\ucDk{D:k:paths}$.}
\begin{proposition}\label{D:p:threats_proba}
    For every $w\in\RR\times\NN$, $k\in\mathcal{S}$ and $r\in\NN^*$, \begin{align}\label{D:estimate_final}
\PP(w\text{ is not $(k,r)$-threatened})\leqslant \ucD{D:c:threats}\,r^{-\alpha}.
\end{align}
\end{proposition}

In order to prove this result, we adapt the proof of \cite[Lemma 5.5]{BHT}, by following again the road map of Section \ref{D:sss:renormalization}. However here it will be a little trickier because we will need to show the subsequence estimate in two steps.

Recall constants defined in \eqref{D:constants_threats}. We are working on the subsequence given by indexes $r=3^j$ for $j\in\NN$. We set $$q_j=q_j^{(k)}=\sup_{w\in \RR\times\NN}\PP(w\text{ is not }(k,3^j)\text{-threatened}).$$
We begin by showing that $q_j$ converges to $0$ when $j\to\infty$, uniformly in $k\in\mathcal{S}$. More precisely, we show the following result.

\begin{lemma}\label{D:estimate2}
    For every $j\geqslant 2$,
\begin{align}
\sup_{k\in\mathcal{S}}\; q_{j_0+j}\leqslant 2^{-j/2}.\label{D:estimate}
\end{align}
\end{lemma}

Note that this bound is not the final bound that we want (which involves $\alpha$). That is why we will need to show a second estimate after this one.

\begin{proof}
In order to show \eqref{D:estimate}, we use induction on $j\geqslant 2$.

\paragraph{Base case.} If a point is not $(k,r)$-threatened, in particular it is not $k$-trapped, so, by \eqref{D:e:for_traps}, we have
    \begin{align*}
    \sup_{k\in\mathcal{S}} \;q_{j_0+2}\leqslant \sup_{k\in\mathcal{S}}\sup_{w\in \R\times\NN} \PP(\text{$w$ is not $k$-trapped})\leqslant 1/2.
    \end{align*}
Therefore, case $j=2$ in \eqref{D:estimate} is satisfied.

\paragraph{Induction step.} Fix $k\in\mathcal{S}$, $j\geqslant 2$ and $w\in \R\times\NN.$ Suppose that \eqref{D:estimate} is satisfied. Note that the event given by $\{\text{$w$ is not $(k,3^{j_0+j+1})$-threatened}\}$ is included in 
both
$$A_j^{(1)}=\bigcap_{i=0}^{3^{j_0+j}-1} \{\text{$w_i$ is not $k$-trapped}\}\;\;\;\text{and}\;\;\;A_j^{(2)}=\bigcap_{i=2\cdot 3^{j_0+j}}^{3^{j_0+j+1}-1} \{\text{$w_i$ is not $k$-trapped}\}.$$ Now, these events are measurable with respect to boxes $B_{3^{j_0+j}h_k}(w)$ and $B_{3^{j_0+j}h_k}(w_{2\cdot 3^{j_0+j}})$ respectively. Let us justify this for the first event. Since $\delta\leqslant 1/2$, for every $w'\in\RR\times\NN$, we have $J_k(w')\subseteq I_{h_k}(w')$. Therefore, using Fact \ref{D:f:allowed_path_box}, $\{\text{$w'$ is $k$-trapped}\}$ is measurable with respect to $B_{h_k}(w')$. Now, for every $i\in\llbracket 0, 3^{j_0+j}-1\rrbracket$, $B_{h_k}(w_i)\subseteq B_{3^{j_0+j}h_k}(w)$, since $|v_+|\leqslant R$ (see \eqref{D:e:bounds_on_v_+}).

Now, $B_{3^{j_0+j}h_k}(w)$ and $B_{3^{j_0+j}h_k}(w_{2\cdot 3^{j_0+j}})$ satisfy the assumptions of Fact \ref{D:p:mixing} with $H=3^{j_0+j}h_k$, so
    \begin{align*}
        q_{j_0+j+1}
        \leqslant q_{j_0+j}^2+\ucD{D:c:mixing}\,h_k^{-\alpha}3^{-(j_0+j)\alpha}
        \leqslant 2^{-j}+\ucD{D:c:mixing}\,3^{-(j_0+j)\alpha},
    \end{align*}
    using the induction assumption and the fact that $h_k\geqslant 1$. In the end, using the fact that $\alpha\geqslant 1$ and \eqref{D:constants_threats},
    \begin{align*}
        \frac{q_{j_0+j+1}}{2^{-(j+1)/2}}
        \leqslant 2^{-(j-1)/2}+\ucD{D:c:mixing}\,2^{(j+1)/2} 3^{-(j_0+j)}
        \leqslant \frac{1}{\sqrt{2}}+\ucD{D:c:mixing} 3^{-j_0}\leqslant 1.
    \end{align*}
    This concludes the induction and therefore the proof of \eqref{D:estimate}.
\end{proof}

We now prove the desired estimate on subsequence $(3^j)_{j>j_1}$; more precisely, we prove the following.

\begin{lemma}
    For every $j\in\NN^*$,
\begin{align}\label{D:pataprout}
\sup_{k\in\mathcal{S}}\;q_{j_1+j}\leqslant \frac{1}{2}\,3^{-\alpha j}.
\end{align}
\end{lemma}

\begin{proof} First note that, using Lemma \ref{D:estimate2} along with \eqref{D:constants_threats}, we have
$$\sup_{k\in\mathcal{S}}\;q_{j_1+1}\leqslant 2^{-(j_1-j_0+1)/2}\leqslant 2^{-(\lceil 2(\alpha\log_2(3)+1)\rceil+1)/2}\leqslant \frac{1}{2}\;3^{-\alpha}.$$

Then, we show \eqref{D:pataprout} by induction on $j\geqslant 1$, using the same method as for the first estimate. For the induction step, we have, using \eqref{D:constants_threats}, \begin{align*}
\frac{q_{j_1+j+1}}{\frac{1}{2}3^{-\alpha(j+1)}}
&\leqslant 2\cdot 3^{\alpha(j+1)}\,\left(\frac{1}{4}3^{-2\alpha j}+\ucD{D:c:mixing} h_k^{-\alpha}3^{-\alpha(j_1+j)}\right)\leqslant \frac{1}{2}3^{-\alpha(j-1)}+2\ucD{D:c:mixing} 3^{-\alpha (j_1-1)}\leqslant 1.
\end{align*}
This concludes the induction and the proof of \eqref{D:pataprout}.
\end{proof}

\begin{proof}[Proof of Proposition \ref{D:p:threats_proba}] We now interpolate \eqref{D:pataprout}. Recall the definition of $\ucD{D:c:threats}$ from \eqref{D:constants_threats}. Let $r\geqslant 3^{j_1+1}$ and $j\in\NN^*$ such that $3^{j_1+j}\leqslant r<3^{j_1+j+1}.$ Then, for every $k\in\mathcal{S}$ and $w\in\RR\times\NN$,
$$\PP(w\text{ is not $(k,r)$-threatened})\leqslant\PP(w\text{ is not $(k,3^{j_1+j})$-threatened})\overset{\eqref{D:pataprout}}{\leqslant} \frac{1}{2} 3^{-\alpha j}\leqslant \frac{3^{\alpha(j_1+1)}}{2} \,r^{-\alpha}\leqslant \ucD{D:c:threats} r^{-\alpha},$$
which proves the result for $r\geqslant 3^{j_1+1}$. Now, for $r<3^{j_1+1}$ and $w\in\RR\times\NN$, we have
$$\PP(w\text{ is not $(k,r)$-threatened})\leqslant 1= \ucD{D:c:threats} 3^{-\alpha(j_1+1)}\leqslant \ucD{D:c:threats} r^{-\alpha}.$$
\end{proof}

\subsection{Threatened density}\label{D:ss:threatened_paths}
We now know that each point has a high probability of being threatened. But what we want is for many points along the sample path of our random walks to be near threatened points. Mind that this is not a direct consequence of Proposition \ref{D:p:threats_proba}, for at first glance it could be the case that our random walk goes precisely to regions where there are not a lot of threatened points.

Let $k\in\NN$. Recall that $l_k=\lfloor L_k^{1/4}\rfloor=L_{k+1}/L_k$ and $h_k=L_{\ucDk{D:k:paths}}L_k$, where $\ucDk{D:k:paths}$ was defined at the beginning of Section \ref{D:ss:traps}. In the rest of the paper, we are going to work with $r=l_{\ucDk{D:k:paths}}$ and the three sets of scales given by $h_k$ (from \eqref{D:hk}), $H_k=rh_k=L_{\ucDk{D:k:paths}+1}L_k$ and $\mathbf{H}_k=L_k^2$. We also set $m_k=\mathbf{H}_k/h_k=L_k/L_{\ucDk{D:k:paths}}$ and $M_k=\mathbf{H}_k/H_k=L_k/L_{\ucDk{D:k:paths}+1}$.

We do not want to look at too many points for entropy reasons, so we will restrict ourselves to rounded points, as defined by the following definition. Recall \eqref{D:d:delta}.

\begin{definition}\label{D:d:rounded_points}
Let $y\in\LLL$ and $k\in\NN$. We denote by $\lfloor y\rfloor_k$ the closest point to the left of $y$ that is in $(\lfloor \delta h_k/4\rfloor\,\ZZ)\times\NN$, which means
$$\lfloor y \rfloor_k=\left(\left\lfloor \frac{\pi_1(y)}{\lfloor \delta h_k/4\rfloor}\right\rfloor \lfloor \delta h_k/4\rfloor,\,\pi_2(y)\right).$$
This is well-defined, because of \eqref{D:condition_on_k0}.
\end{definition}

\begin{definition}\label{D:d:density}
Let $k\in\mathcal{S}$ and $(\chi(n))_{n\in\llbracket 0,\mathbf{H}_k\rrbracket}$ an allowed path starting in $I_{\mathbf{H}_k}$. We define its threatened density as 
\begin{align}\label{D:e:density}
D_k(\chi)=\frac{1}{M_k} \,\#\left\{0\leqslant j<M_k,\;\lfloor \chi(jH_k)\rfloor_k \text{ is $(k,r)$-threatened}\right\}.\end{align}
\end{definition}

\ncD{D:c:density}
\begin{proposition}\label{D:p:density}
There exists $\ucD{D:c:density}>0$ such that for every $k\in\mathcal{S}$ satisfying $k>\ucDk{D:k:paths}$, we have
$$\PP\left(\exists\,(\chi(n))_{n\in\llbracket 0,\mathbf{H}_k\rrbracket}\;\text{allowed path starting in $I_{\mathbf{H}_k}$},\,D_k(\chi)<1/2\right)\leqslant \ucD{D:c:density}L_k^{-(\alpha-1)/5}.$$
\end{proposition}

\begin{proof}
The trick is that in $\mathbf{H}_k=L_k\cdot L_k$, we are going to fix the first $k$ {\color{black}(it is noted $\hat{k}$ from now on)}, and use an induction on the second $k$, in order to show a stronger result. For every $w\in\RR\times\NN$, if $(\chi(n))_{n\in\llbracket 0,L_{\hat{k}}L_k\rrbracket}$ is an allowed path starting in $I_{L_{\hat{k}}L_k}(w)$, we set $$D_{\hat{k},k}(\chi)=\frac{1}{M_k} \,\#\left\{0\leqslant j<M_k,\;\lfloor \chi(jH_{\hat{k}})\rfloor_{\hat{k}} \text{ is $(\hat{k},r)$-threatened}\right\}.$$
Note that if $w=(0,0)$ and $k=\hat{k}$, then $D_k(\chi)=D_{\hat{k},k}(\chi).$ We aim at showing that there exists $\ucD{D:c:density}>0$ such that for every $\hat{k}\in\mathcal{S}$ and $k>\ucDk{D:k:paths}$, we have
$$\PP\left(\exists\,(\chi(n))_{n\in\llbracket 0,L_{\hat{k}}L_k\rrbracket}\;\text{allowed path starting in $I_{L_{\hat{k}}L_k}$},\,D_{\hat{k},k}(\chi)<1/2\right)\leqslant \ucD{D:c:density}L_k^{-(\alpha-1)/5}.$$
Now, the rest of the proof is very similar to the proof of Lemma \ref{D:l:bounds_on_p} (see Section \ref{D:ss:deviation_proof}): we are going to use the renormalization method 
from Section \ref{D:sss:renormalization}, along with a sequence of densities $(\rho_k)_{k\geqslant \ucDk{D:k:paths}}$. We first fix $\hat{k}\in\mathcal{S}$, and we define 
\begin{align}\label{D:d:sequence_densities}
    \left\{\begin{array}{l}
    \rho_{\ucDk{D:k:paths}}=1\\
    \forall k\geqslant \ucDk{D:k:paths},\;\;\rho_{k+1}=\rho_k-\frac{2}{l_k}.
    \end{array}\right.
\end{align}
We can check that by definitions of the $(l_k)$ in \eqref{D:e:scalesL_k}, we have $\sum_{k\geqslant 1} \frac{2}{l_k}\leqslant \frac{1}{2}.$ Therefore, $\rho_k\geqslant 1/2$ for every $k\geqslant \ucDk{D:k:paths}.$ We define, for $w\in\RR\times\NN$ and $k\in\NN$,
$$S_k(w)=S_k^{\hat{k}}(w)=\left\{\exists\,(\chi(n))_{n\in\llbracket 0,L_{\hat{k}}L_k\rrbracket}\text{ allowed path starting in }I_{L_{\hat{k}}L_k}(w),\;D_{\hat{k},k}(\chi)\leqslant\rho_k\right\}$$
and we set $S_k=S_k((0,0))$ as usual. Since $\rho_k\geqslant 1/2$, it will be sufficient to show that for every $k>\ucDk{D:k:paths}$,
\begin{align}\label{D:e:estimate_densities}
    \PP(S_k)\leqslant \ucD{D:c:threats} \frac{4^{\alpha+2}}{\delta}\, L_k^{-(\alpha-1)/5}.
\end{align}

To do this, we use induction on $k>\ucDk{D:k:paths}$.
\paragraph{Base case.} Using Proposition \ref{D:p:threats_proba} as well as the fact that $M_{\ucDk{D:k:paths}+1}=1$ and $\rho_{\ucDk{D:k:paths}+1}<1$, we have
    \begin{align*}
    \PP\left(S_{\ucDk{D:k:paths}+1}\right)
    &=\PP\left(\exists\,y\in I_{L_{\hat{k}}L_{\ucDk{D:k:paths}+1}},\,\lfloor y\rfloor_{\hat{k}}\text{ is not $(\hat{k},r)$-threatened}\right)\\
    &\leqslant \left\lceil\displaystyle\frac{L_{\hat{k}}L_{\ucDk{D:k:paths}+1}}{\left\lfloor \frac{\delta L_{\hat{k}}L_{\ucDk{D:k:paths}}}{4}\right\rfloor}\right\rceil \ucD{D:c:threats}\,r^{-\alpha}\leqslant \ucD{D:c:threats} \frac{4^{\alpha+2}}{\delta}\,L_{\ucDk{D:k:paths}+1}^{-(\alpha-1)/5}
    \end{align*}
(in the last bound, we used very rough inequalities for lower and upper integral parts, for instance $l_{\ucDk{D:k:paths}}\geqslant \frac{1}{2}L_{\ucDk{D:k:paths}}^{1/4}$). Therefore, we do get inequality \eqref{D:e:estimate_densities} with $k=\ucDk{D:k:paths}+1$. Note that considering only rounded points $\lfloor y\rfloor_{\hat{k}}$ was crucial here to obtain a bound that is uniform in $\hat{k}$.

\paragraph{Induction step.} We now fix $k>\ucDk{D:k:paths}$ and assume that \eqref{D:e:estimate_densities} is true. We want to show that the same estimate holds for $k+1$. Recall the definition of $\CCC$ from \eqref{D:CCC}: here we work with $\CCC=\CCC_{L_{\hat{k}},k}$. We claim that on $S_{k+1}$, there exist $w_1,\,w_2\in\CCC$ such that
    \begin{align}\label{D:e:S_cascade}\left\{\begin{array}{l} S_{k}(w_1)\text{ and }S_{k}(w_2)\text{ occur};\\
    \vert\pi_2(w_1)-\pi_2(w_2)\vert\geqslant 2L_{\hat{k}}L_k.\end{array}\right.\end{align}
Indeed, assume that $S_{k+1}$ occurs. This means that there exists an allowed path $(\chi(n))_{n\in\llbracket 0,L_{\hat{k}}L_{k+1}\rrbracket}$ starting in $I_{L_{\hat{k}}L_{k+1}}$ such that $D_{\hat{k},k+1}(\chi)\leqslant\rho_k$. Now, assume by contradiction that we cannot find three points $w_1',\,w_2',\,w_3'\in \CCC$ such that the $(\pi_2(w_i'))_{i=1,2,3}$ are distinct and $S_{k}(w_i')$ occurs for every $i\in\{1,2,3\}$. This implies that for $l_k-2$ values of $i\in\llbracket 0,l_k-1\rrbracket$, we have 
$$\#\left\{0\leqslant j<M_k,\,\lfloor \chi(iL_{\hat{k}}L_k+jH_{\hat{k}})\rfloor_{\hat{k}}\text{ is $(\hat{k},r)$-threatened}\right\}>\rho_k M_k.$$
Therefore, recalling Fact \ref{D:f:allowed_path_box},
    \begin{align*}
        D_{\hat{k},k+1}(\chi)
        &=\frac{1}{M_{k+1}} \,\#\left\{0\leqslant j<M_{k+1},\;\lfloor \chi(jH_{\hat{k}})\rfloor_{\hat{k}} \text{ is $(\hat{k},r)$-threatened}\right\}\\
        &>\frac{1}{M_{k+1}}\,(l_k-2)M_k\rho_k=\frac{l_k-2}{l_k}\rho_k> \rho_{k+1},
    \end{align*}
    which contradicts the fact that $S_{k+1}$ occurs. Therefore, on $S_{k+1}$, \eqref{D:e:S_cascade} holds for two points $w_1,w_2\in\mathcal{C}$. Now, for $i\in\{1,2\}$, $S_{k}(w_i)$ is measurable with respect to box $B_{L_{\hat{k}}L_k}(w_i)$, and the two boxes satisfy the assumptions of Fact \ref{D:p:mixing} with $H=L_{\hat{k}}L_k$. Therefore, using the fact that $\PP(S_k(w))=\PP(S_k)$ for any $w\in\RR\times\NN$, we have
    $$\PP(S_{k+1})\leqslant |\CCC|^2 (\PP(S_k)^2+\ucD{D:c:mixing} L_k^{-\alpha})=(2R+3)^2\, l_k^4\, (\PP(S_k)^2+\ucD{D:c:mixing} L_k^{-\alpha}).$$
    In the end, using the induction assumption,
    \begin{align*}
        \frac{\PP(S_{k+1})}{L_{k+1}^{-(\alpha-1)/5}}
        &\leqslant L_k^{(\alpha-1)/4} (2R+3)^2\, l_k^4\, (\PP(S_{k})^2+\ucD{D:c:mixing}L_k^{-\alpha})\\
        &\leqslant (2R+3)^2\, L_k^{(\alpha+3)/4}\,\left(\ucD{D:c:threats}^2 \frac{4^{2\alpha+4}}{\delta^2} L_k^{-2(\alpha-1)/5} +\ucD{D:c:mixing}\,L_k^{-\alpha}\right)\\
        &\leqslant (2R+3)^2\,\left(\ucD{D:c:threats}^2\frac{4^{2\alpha+4}}{\delta^2}+\ucD{D:c:mixing}\right) L_k^{-(3\alpha-23)/20}\\
        &\leqslant \ucD{D:c:threats} \frac{4^{\alpha+2}}{\delta}&\mbox{using \eqref{D:e:condition_on_k1}.}
    \end{align*}
    This concludes the induction. Now, by setting $\ucD{D:c:density}=\ucD{D:c:threats}\frac{4^{\alpha+2}}{\delta}$ and taking $k\in\mathcal{S}$ as well as $\hat{k}=k$, the proposition is shown.
\end{proof}

\subsection{Delays}\label{D:ss:delays}

Keep in mind that our goal is to prove \eqref{D:e:contradiction}. Recall \eqref{D:hk} and the definitions from the beginning of Section \ref{D:ss:threatened_paths}. For $k\in\mathcal{S}$ satisfying $k>\ucDk{D:k:paths}$ and $y\in I_{\mathbf{H}_k}$, the points we are going to consider along the sample path of $Z^y$ are the following.
\begin{itemize}
    \item Potential threats:
    \begin{align}\label{D:d:W_j}
    W_j^y=\lfloor Z^y_{jH_k}\rfloor_k\;\;\;\text{for}\;j\in\NN;
    \end{align}
    \item Potential traps:
    \begin{align}\label{D:d:w_j}
w_j^y=W_{j_0}^y+j_1h_k(v_+,1)\;\;\;\text{for}\;j=j_0r+j_1\in\NN,\text{ with } j_0\in\NN \text{ and } 0\leqslant j_1<r.
\end{align}
\end{itemize}

We let $(\Tilde{U}_j)_{j\in\NN}=(\Tilde{U}_j(k))_{j\in\NN}$ be i.i.d. uniform random variables in $\llbracket 0,\lfloor \delta h_k\rfloor-1\rrbracket\times\{0\}$, which are independent of all the other random variables introduced so far. For every $j\in\NN$, we set $y_j=y_j(y,k)\in\LLL$ to be the leftmost point of $J_k(w_j^y)$ (recall definitions \eqref{D:d:N_k}). We set
\begin{align}\label{D:d:Y_j}
Y_j=Y_j(y,k)=y_j+\Tilde{U}_j.
\end{align}
Therefore $y_j+\llbracket 0,\lfloor \delta h_k\rfloor-1\rrbracket\times\{0\}$ is included in $J_k(w_j^y)$ but might exclude its rightmost point.

We keep notation $\PP$ for the probability on the larger space that we created by considering these new random variables. We let $$\PP_{c}=\PP\left(\cdot\,\vert\,\omega,\,(U_n^x)_{n\in\N,x\in\ZZ}\right).$$

Let $j\in\NN.$ Recall the definition of $v_0$ from \eqref{D:d:delta} and $m_k$ from the beginning of Section \ref{D:ss:threatened_paths}. We define the following events and corresponding sets:
\begin{align*}\begin{array}{ll}
    \mathrm{Trap}_j=\left\{w_j^y\text{ is $k$-trapped}\right\}& \mathrm{Trap}=\{0\leqslant j<m_k,\,\mathrm{Trap}_j\text{ occurs}\};\\
    \mathrm{Bar}_j=\left\{V_{h_k}^{Y_j}\leqslant v_0\right\}& \mathrm{Bar}=\{0\leqslant j<m_k,\,\mathrm{Bar}_j\text{ occurs}\};\\
    \mathrm{Cross}_j=\left\{X^y_{jh_k}<\pi_1(Y_j),\,X_{(j+1)h_k}^y>X_{h_k}^{Y_j}\right\}&\mathrm{Cross}=\{0\leqslant j<m_k,\,\mathrm{Cross}_j\text{ occurs}\};\\
    \mathrm{Del}_j=\left\{X_{(j+1)h_k}^y\leqslant \pi_1(w_j^y)+(v_+-\delta)h_k\right\}& \mathrm{Del}=\{0\leqslant j<m_k,\,\mathrm{Del}_j\text{ occurs}\}.\\
    \end{array}
\end{align*}
Here we chose not to write the dependencies on $k$ and $y$ for simplicity, but when it makes things clearer, we will write for instance $\mathrm{Del}=\mathrm{Del}(y,k)$.

Let us give some motivation for these definitions, illustrated by Figure \ref{D:f:density}. $\mathrm{Bar}_j$ means that the sample path starting at the random point $Y_j$ may act as a barrier for our random walk $Z^y$ between times $jh_k$ and $(j+1)h_k$. It is possible to cross this barrier, which is roughly what event $\mathrm{Cross}_j$ stands for. But it may also be the case that because of this barrier, our random walk is delayed at scale $h_k$, which is what $\mathrm{Del}_j$ stands for.

\begin{figure}[!h]
  \begin{center}
    \begin{tikzpicture}[use Hobby shortcut, scale=0.7]
        \draw[thick] (0,0) -- (5,0);
        \draw[thick] (0,0.1) -- (0,-0.1);
        \draw[thick, below left] (5,0.1) -- (5,-0.1) node {$I_{\mathbf{H}_k}$};
        \draw[thick, below left] (-5,5) rectangle (10,0);
        \draw[below left] (9,4.8) node {$B_{\mathbf{H}_k}$};
        \draw[fill, above right] (2,0) circle (.05) node {$y$};
        \draw (2,0) .. (1.5,1) .. (2.5,2) ..(2.1,2.7) .. (2,2.8).. (2,3) .. (3.3,4) .. (3,5);
        
        \draw[dashed] (-5,1) -- (10,1);
        \draw[dashed] (-5,2) -- (10,2);
        \draw[dashed] (-5,3) -- (10,3);
        \draw[dashed] (-5,4) -- (10,4);
        \draw[<->,>=latex] (-4,2) -- (-4,3);
        \draw[right] (-4,2.5) node {$H_k$};
        
        \draw[fill, below left] (1.7,0) circle (.1) node{$\lfloor y\rfloor_k=W_0^y$};
        \draw[fill=white, above left] (1.4,1) circle (.1) node {$W^y_1$};
        \draw[fill, below right] (2.5,2) circle (.1) node {$W^y_2$};
        \draw[fill=white, below left] (1.8,3) circle (.1);
        \draw[fill, below left] (3.1,4) circle (.1);
        
        \draw[densely dotted,very thick] (2.6,2.5) .. (2.4,2.6) .. (2.1,2.7) .. (2,2.8);
        \draw[fill, right] (2.6,2.5) circle (.04) node{$Y_j$};
        \draw[fill] (2,2.8) circle (.04);
        \draw[<->] (1.6,2.45) -- (1.6,2.85) ;
        \draw[left] (1.6,2.65) node {$h_k$};
    \end{tikzpicture}
    \caption{Along the way, the random walk starting at $y$ passes near threatened points represented by the big filled dots. These points are likely to be responsible for a delay of our random walk. Here, the second threatened point encountered comes with a trap that causes a delay. Events $\mathrm{Trap}_{j}$, $\mathrm{Bar}_{j}$ and $\mathrm{Del}_j$ occur (for some $j\in\llbracket 2r,3r-1\rrbracket$), in fact here our random walk coalesces with $Z^{Y_j}$.}
    \label{D:f:density}
    \end{center}
\end{figure}
Our goal is to show a lower bound on the cardinality of $\mathrm{Del}$ with high probability. More precisely, we aim at showing that a positive proportion of indexes $j$ satisfy $\mathrm{Del}_j$, in order to delay our random walks. This is shown in Lemma \ref{D:p:delay}. In order to get there, we will use several other lemmas that deal with the other events introduced above.

The first lemma gives a lower bound on $|\mathrm{Bar}|$ with high probability.

\ncD{D:c:bar}
\begin{lemma}\label{D:l:bar}
There exists $\ucD{D:c:bar}>0$ such that for every $k\in\mathcal{S}$ satisfying $k>\ucDk{D:k:paths}$ and $y\in I_{\mathbf{H}_k}$,
$$\PP\left(|\mathrm{Bar}(y,k)|<\frac{M_k}{10}\right)\leqslant \ucD{D:c:bar} L_k^{-(\alpha-1)/5}.$$
\end{lemma}

\begin{proof}
    We first write    \begin{align}\label{D:e:2_parts}\PP\left(|\mathrm{Bar}|<\frac{M_k}{10}\right)\leqslant \PP\left(|\mathrm{Trap}|<\frac{M_k}{2}\right)+\PP\left(|\mathrm{Trap}|\geqslant \frac{M_k}{2},\;|\mathrm{Bar}|<\frac{M_k}{10}\right).\end{align}
To control the first term on the right-hand side of \eqref{D:e:2_parts}, note that the number of threatened points among $\{W_j^y,0\leqslant j< M_k\}$ is at most the number of traps among $\{w_j^y,0\leqslant j< m_k\}$ (recall \eqref{D:d:W_j}, \eqref{D:d:w_j} and Definition \ref{D:d:threats}). Therefore, recalling Definition \ref{D:d:density}, we have
\begin{align*}
        \PP\left(|\mathrm{Trap}|<\frac{M_k}{2}\right)
        &\leqslant \PP\left(D_k\left((Z^y_n)_{0\leqslant n\leqslant\mathbf{H}_k}\right)<1/2\right)\\
        &\leqslant \PP\left(\exists\,(\chi(n))_{n\in\llbracket 0,\mathbf{H}_k\rrbracket}\;\text{allowed path starting in $I_{\mathbf{H}_k}$},\,D_k(\chi)<1/2\right)\\
        &\leqslant \ucD{D:c:density}L_k^{-(\alpha-1)/5},
    \end{align*}
    using Proposition \ref{D:p:density}.
    As for the second term in \eqref{D:e:2_parts}, it is at most
    \begin{align*}
        \PP\left(|\mathrm{Bar}|<\frac{|\mathrm{Trap}|}{5},\,|\mathrm{Trap}|\geqslant \frac{M_k}{2}\right)
        &\leqslant \sum_{J\subseteq \llbracket 0,m_k-1\rrbracket\atop |J|\geqslant M_k/2} \PP\left(\mathrm{Trap}=J,\;|\mathrm{Bar}|<\frac{|J|}{5}\right)\\
        &=\sum_{J\subseteq \llbracket 0,m_k-1\rrbracket\atop|J|\geqslant M_k/2} \EE\left[\mathbf{1}_{\mathrm{Trap}=J}\;\PP_c\left(|\mathrm{Bar}|<\frac{|J|}{5}\right)\right]
    \end{align*}
    Now, using Definition \ref{D:d:traps}, note that for every $J$, on $\{\mathrm{Trap}=J\}$, we have, for every $j\in J$,
    \begin{align}
\PP_c(\mathrm{Bar}_j)
&\geqslant \frac{1}{\lfloor \delta h_k\rfloor}\, \#\{z\in y_j+\llbracket 0,\lfloor \delta h_k\rfloor-1\rrbracket\times\{0\},\,V^z_{h_k}\leqslant v_0\}\nonumber\\
&\geqslant \frac{N_k(w_j^y)/3-1}{\lfloor \delta h_k\rfloor} \nonumber=\frac{1}{3}-\frac{1}{\lfloor\delta h_k\rfloor}\nonumber\\
&\geqslant 1/4.\label{D:ineq}
    \end{align}
    In the last inequality, we used \eqref{D:condition_on_k0}. Now, note that under $\PP_c$, events $\mathrm{Bar}_j$ are independent, because the $(\tilde{U}_j)$ are. Therefore, if $a>0$ is fixed, we have, for any $J\subseteq \llbracket 0,m_k-1\rrbracket$ such that $|J|\geqslant M_k/2$ and on $\{\mathrm{Trap}=J\}$,
    \begin{align*}
    \PP_c\left(|\mathrm{Bar}|<\frac{|J|}{5}\right)
    &\leqslant \PP_c\left(\sum_{j\in J} \mathbf{1}_{\mathrm{Bar}_j^c}\geqslant \frac{4|J|}{5}\right)&\mbox{using that $\mathrm{Bar}\supseteq \mathrm{Bar}\cap J$}\\
    &=\PP_c\left(\prod_{j\in J} \exp(a\mathbf{1}_{\mathrm{Bar}_j^c})\geqslant \exp(4a|J|/5)\right)\\
    &\leqslant e^{-4a|J|/5} \left(\frac{3}{4}(e^a-1)+1\right)^{|J|}&\mbox{using \eqref{D:ineq} and independence}\\
    &=\left(1-\frac{a}{20}+o_{a\to 0}(a)\right)^{|J|}\\
    &\leqslant e^{-c|J|}&\mbox{for a good choice of $a$}\\
    &\leqslant e^{-cL_k}&\mbox{since $|J|\geqslant M_k/2=cL_k$.}
    \end{align*}
    Hence the result by choosing $\ucD{D:c:bar}$ properly, independently of $y$ and $k$.
\end{proof}

The second lemma that we are going to show roughly gives a uniform lower bound on the probability of the crossing event. This lemma is crucial, since its goal is to replace the monotonicity property from \cite{BHT} that we are lacking when $R\geqslant 2$ (see Remark \ref{D:r:monotonicity_BHT}).

Recall definition \eqref{D:d:Y_j}. We consider the filtration given for $j\in\NN$ by $$\mathcal{F}_j=\mathcal{F}_j(y,k)=\sigma\left(\{Z^y_n,\,0\leqslant n\leqslant jh_k\}\cup\{Z^{Y_i}_n,\,0\leqslant i\leqslant j,\,0\leqslant n\leqslant h_k\}\right).$$

It will also be interesting to consider all the indexes $j$ such that $\mathrm{Bar}_j$ occurs, without stopping at $m_k-1$. We define a sequence of $(\mathcal{F}_j)$-stopping times $(\tau_i)_{i\in\NN}$ by setting $\tau_{-1}=-1$ and for all $i\in\NN$,
\begin{align*}
    \tau_i=\left\{\begin{array}{ll}\inf\{j>\tau_{i-1},\,\mathrm{Bar}_j\text{ occurs}\}&\text{if this set in non-empty}\\ +\infty&\text{otherwise}\end{array}\right.
\end{align*}
We consider the stopped filtrations $(\mathcal{F}_{\tau_i})_{i\in\NN}$, and we set $\mathrm{Cross}_\infty=\emptyset.$

Recall constant $\gamma>0$ from Assumption \ref{D:a:uniform_ellipticity}.

\begin{lemma}\label{D:l:unif_bound_gamma}
For every $i\in\NN$, {\color{black}almost surely} we have
\begin{align}\label{D:e:lemma2}
    \PP(\mathrm{Cross}_{\tau_i}\vert \mathcal{F}_{\tau_i})\leqslant 1-\gamma.
\end{align}
\end{lemma}

\begin{proof}
Let us start by showing that for every $j\in\NN$,
\begin{align}\label{D:e:lemma1}
{\color{black}\text{a.s.}}\quad \PP(\mathrm{Cross}_j\,\vert\,\mathcal{F}_j)\leqslant 1-\gamma.
\end{align}

{\color{black}
Let $j\in\N$. Here, the idea is to consider several scenarios in which $\mathrm{Cross}_j$ cannot occur. Two good scenarios are when $Z^y$ stays too far from $Z^{Y_j}$ to have any chance of crossing its path, and when the two paths coincide from the start ($Z^y_{jh_k}=Y_j$). The bad scenario is when $Z^y$ approaches $Z^{Y_j}$ close enough to have a chance of crossing its path, but then the two paths also have a chance to coalesce.

Let us define
$$T=\left\{\begin{array}{ll}\min\left\{n\geqslant jh_k,\,\left\vert X^y_{n}-X^{Y_j}_{n+1-jh_k}\right\vert \leqslant R\right\}&\text{if this set in non-empty;}\\ +\infty&\text{otherwise,}\end{array}\right.$$
and the following event:
$$A=\left\{g\left(\omega_{T}(X_{T}^y-\ell),\ldots,\omega_{T}(X_{T}^y+\ell),U_{T}^{X_{T}^y}\right)=X_{T+1-jh_k}^{Y_j}-X_{T}^y\right\}.$$
Suppose that we have $T<(j+1)h_k$ and that $A$ occurs; in other words it is possible for $Z^y$ to jump to location $Z^{Y_j}_{T+1-jh_k}$ between times $T$ and $T+1$, since $\vert X^y_{T}-X^{Y_j}_{T+1-jh_k}\vert\leqslant R$. If this happens, then both paths coalesce, because of \eqref{D:e:coupling_property}, so $\mathrm{Cross}_j$ does not occur. Furthermore, when $T\geqslant (j+1)h_k$ and $Z^y_{jh_k}=Y_j$, $\mathrm{Cross}_j$ cannot occur either.

Let $\zeta$ be a non-negative $\mathcal{F}_{j}$-measurable random variable. Using the three incompatible scenarios discussed above, we get
\begin{align}\label{D:last_term}
    \EE\left[\mathbf{1}_{\mathrm{Cross}_j^c} \zeta\right]
    &\geqslant \EE\left[\mathbf{1}_{T<(j+1)h_k}\, \mathbf{1}_{A}\,\mathbf{1}_{Z^y_{jh_k}\neq Y_j} \zeta \right]+\EE\left[\mathbf{1}_{T\geqslant (j+1)h_k}\zeta\right]+\EE\left[\mathbf{1}_{Z^y_{jh_k}=Y_j}\zeta\right].
\end{align}
Note that we put aside the $Z^y_{jh_k}=Y_j$ case in the first term, for in that case the argument of independence we are about to use falls apart. The first term on the right-hand side in \eqref{D:last_term} can be rewritten as
\begin{align*}
    \displaystyle \sum_{\substack{|x_2-x_1|\leqslant R\\\sigma_{-\ell},\ldots,\sigma_{\ell}\in S\\ jh_k\leqslant t<(j+1)h_k}} \EE\left[\mathbf{1}_{g\left(\sigma_{-\ell},\ldots,\sigma_{\ell},U_t^{x_1}\right)=x_2-x_1}\, \mathbf{1}_{T=t}\,\mathbf{1}_{\omega_t(x_1-\ell)=\sigma_{-\ell},\ldots,\atop \omega_t(x_1+\ell)=\sigma_\ell} \,\mathbf{1}_{X^y_t=x_1,\,X^{Y_j}_{t+1-jh_k}=x_2}\,\mathbf{1}_{Z^y_{jh_k}\neq Y_j}\,\zeta \right].
\end{align*}
In each term of this sum, $U_t^{x_1}$ is independent of $$\mathbf{1}_{T=t}\,\mathbf{1}_{\omega_t(x_1-\ell)=\sigma_{-\ell},\ldots,\omega_t(x_1+\ell)=\sigma_\ell} \,\mathbf{1}_{X^y_t=x_1,\,X^{Y_j}_{t+1-jh_k}=x_2}\,\mathbf{1}_{Z^y_{jh_k}\neq Y_j}\,\zeta.$$ More precisely, the latter random variable is measurable with respect to $\omega$ and the uniform variables $U_n^x$ with $n\neq t$ or $x\neq x_1$. Indeed, assuming that $T=t$, $X_t^y=x_1$ and $Z^y_{jh_k}\neq Y_j$, $Z^{Y_j}$ cannot visit $(x_1,t)$, otherwise we would have $X^{Y_j}_{t-jh_k}=x_1=X_t^y$, and so $|X_{t-1}^y-X^{Y_j}_{t-jh_k}|=|X_{t-1}^y-X_t^y|\leqslant R$. Moreover, $t>jh_k$, since $Z^y_{jh_k}\neq Y_j$, so $T\leqslant t-1$, which is a contradiction.

Therefore, using Assumption \ref{D:a:uniform_ellipticity}, we can rewrite the first term in the right-hand site of \eqref{D:last_term} as
\begin{align*}
    \displaystyle \sum_{\substack{|x_2-x_1|\leqslant R\\\sigma_{-\ell},\ldots,\sigma_{\ell}\in S\\ jh_k\leqslant t<(j+1)h_k}}
    \PP\left(g\left(\sigma_{-\ell},\ldots,\sigma_{\ell},U_t^{x_1}\right)=x_2-x_1\right)  \,\EE\left[\mathbf{1}_{T=t}\,\mathbf{1}_{\omega_t(x_1-\ell)=\sigma_{-\ell},\ldots,\atop \omega_t(x_1+\ell)=\sigma_\ell} \,\mathbf{1}_{X^y_t=x_1,X^{Y_j}_{t+1-jh_k}=x_2}\,\mathbf{1}_{Z^y_{jh_k}\neq Y_j}\,\zeta \right]
\end{align*}
which is at least $\gamma \,\EE\left[\mathbf{1}_{T<(j+1)h_k}\,\mathbf{1}_{Z^y_{jh_k}\neq Y_j}\,\zeta\,\right].$\\
In the end, putting together the three terms in \eqref{D:last_term}, we get $\EE\left[\mathbf{1}_{\mathrm{Cross}_j^c}\, \zeta\right]\geqslant \gamma\, \EE[\zeta],$ hence \eqref{D:e:lemma1}.
}

In order to show \eqref{D:e:lemma2}, we fix $i\in\NN$ and we let $\zeta$ be a non-negative $\mathcal{F}_{\tau_i}$-measurable random variable. We have
\begin{align*}
    \EE\left[\mathbf{1}_{\mathrm{Cross}_{\tau_i}}\zeta\right]
    &=\sum_{j\in\NN} \EE\left[\mathbf{1}_{\mathrm{Cross}_j}\mathbf{1}_{\tau_i=j}\zeta\right]&\mbox{because $\mathrm{Cross}_\infty=\emptyset$}\\
    &=\sum_{j\in\NN}\EE\left[\PP\left(\mathrm{Cross}_j\vert \mathcal{F}_j\right)\mathbf{1}_{\tau_i=j}\zeta\right]&\mbox{since $\mathbf{1}_{\tau_i=j}\zeta$ is $\mathcal{F}_j$-measurable}\\
    &\leqslant (1-\gamma)\,\EE[\zeta]&\mbox{using \eqref{D:e:lemma1}.}
\end{align*}
This concludes the proof of \eqref{D:e:lemma2}.
\end{proof}

The third and last lemma that we will need links the events $\mathrm{Bar}_j$ and $\mathrm{Cross}_j$ that we studied in the first two lemmas, to the event $\mathrm{Del}_j$ that we are interested in. However, we first need to define a new event $\mathcal{G}$ of high probability that guarantees that random walks do not go too fast to the right, so that they have a chance of getting delayed by traps. Let \begin{align}\label{D:d:lambda}\lambda=\frac{\gamma\delta}{40 r^2}.\end{align}
Note that
\begin{align}\label{D:e:lam}
\lambda\leqslant\frac{\delta}{2r},
\end{align}
since $\gamma\leqslant 1$ and $r=l_{\ucDk{D:k:paths}}\geqslant 1$. Similarly to \eqref{D:CCC}, we define a set
$$\hat{\CCC}=\hat{\CCC}_k=\{(ih_k,jh_k)\in h_k\ZZ^2,\,-(R+1)m_k\leqslant i<(R+2)m_k,\,0\leqslant j<m_k\},$$
which has cardinality $(2R+3)m_k^2$ and satisfies
\begin{align}\label{D:pocpoc}\underset{w\in\hat{\CCC}}{\bigcup} I_{h_k}(w)=B_{\mathbf{H}_k}\cap(\R\times h_k \Z).\end{align}
Recall definitions \eqref{D:d:Ap}. We set
\begin{align}\label{D:G}\mathcal{G}=\mathcal{G}_k=\bigcap_{w\in \hat{\CCC}} A_{h_k,w}(v_++\lambda)^c.\end{align}
Using Lemma \ref{D:l:bounds_on_p}, we have
\begin{align*}
    \PP(\mathcal{G}^c)
    \leqslant |\hat{C}|\;p_{h_k}(v_++\lambda)
    \leqslant (2R+3)\,m_k^2\,\ucD{D:c:deviation}(\lambda)\, h_k^{-\alpha/4}
    \leqslant  c\left(\frac{L_k}{L_{\ucDk{D:k:paths}}}\right)^2\,(L_k L_{\ucDk{D:k:paths}})^{-\alpha/4}.
\end{align*}
Therefore, using that $\alpha>8,$ \begin{align}\label{D:e:proba_G}
    \PP(\mathcal{G}_k^c)\xrightarrow[k\to\infty]{}0.
\end{align}

\begin{lemma}\label{D:l:useful}
{\color{black}For every $j\in\N$,} we have $$\mathcal{G}\cap\mathrm{Bar}_j\cap\mathrm{Cross}_j^c\subseteq \mathrm{Del}_j.$$
\end{lemma}

\begin{proof}
We are first going to show that
\begin{align}\label{D:e:useful1}
\mathcal{G}\subseteq \bigcap_{j\in \llbracket 0,m_k-1\rrbracket} \left\{X^y_{jh_k}<\pi_1(Y_j)\right\}.
\end{align}
Let $j\in \llbracket 0,m_k-1\rrbracket$, $j=j_0r+j_1$ where $0\leqslant j_0<M_k$ and $0\leqslant j_1<r$ (note that $rM_k=m_k$). On $\mathcal{G}$, we have
\begin{align*}
    X^y_{jh_k}=X_{j_1h_k}^{Z^y_{j_0H_k}}
    &\leqslant X^y_{j_0H_k}+j_1 h_k(v_++\lambda)&\mbox{by \eqref{D:pocpoc}, \eqref{D:G} and Fact \ref{D:f:allowed_path_box}}\\
    &\leqslant X^y_{j_0H_k}+j_1 h_k\left(v_++\frac{\delta}{2r}\right)&\mbox{using \eqref{D:e:lam}}\\
    &\leqslant \pi_1(W_{j_0}^y)+\frac{\delta h_k}{4}+j_1 h_k\left(v_++\frac{\delta}{2r}\right)&\mbox{by \eqref{D:d:W_j} and Definition \ref{D:d:rounded_points}}\\
    &<\pi_1(W_{j_0}^y)+j_1 v_+ h_k +\delta h_k\\
    &=\pi_1(w_j^y)+\delta h_k\\
    &\leqslant \pi_1(y_j)\leqslant \pi_1(Y_j).
\end{align*}
This proves \eqref{D:e:useful1}. Therefore, for a fixed $j\in\llbracket 0,m_k-1\rrbracket$, by definition of $\mathrm{Cross}_j$, we have \begin{align}\label{D:e:useful11}
\mathcal{G}\cap \mathrm{Cross}_j^c\subseteq\left\{X^y_{(j+1)h_k}\leqslant X^{Y_j}_{h_k}\right\}.
\end{align}
This implies that on $\mathcal{G}\cap\mathrm{Bar}_j\cap\mathrm{Cross}_j^c$, we have \begin{align*}
    X^y_{(j+1)h_k}
    &\leqslant X_{h_k}^{Y_j}&\mbox{using \eqref{D:e:useful11}}\\
    &\leqslant \pi_1(Y_j)+v_0h_k&\mbox{by definition of $\mathrm{Bar}_j$}\\
    &\leqslant \pi_1(w_j^y)+2\delta h_k+v_0 h_k&\mbox{by \eqref{D:d:Y_j}}\\
    &\leqslant \pi_1(w_j^y)+(v_+-\delta)h_k&\mbox{by \eqref{D:d:delta},}
\end{align*}
so we are on $\mathrm{Del}_j$.
\end{proof}

\ncD{D:c:del}
\begin{lemma}\label{D:p:delay}
There exists $\ucD{D:c:del}>0$ such that for every $k\in\mathcal{S}$ such that $k>\ucDk{D:k:paths}$ and $y\in I_{\mathbf{H}_k}$, 
$$\PP\left(\mathcal{G}_k,\,|\mathrm{Del}(y,k)|<\frac{\gamma M_k}{20}\right)\leqslant \ucD{D:c:del} L_k^{-(\alpha-1)/5}.$$
\end{lemma}

\begin{proof}
We first write
$$\PP\left(\mathcal{G},\,|\mathrm{Del}|<\frac{\gamma M_k}{20}\right)\leqslant \PP\left(|\mathrm{Bar}|<\frac{M_k}{10}\right)+\PP\left(\mathcal{G},\,|\mathrm{Bar}|\geqslant \frac{M_k}{10},\;|\mathrm{Del}|<\frac{\gamma}{2}\,|\mathrm{Bar}|\right).$$
Using Lemma \ref{D:l:bar}, the first term on the right-hand side of the above inequality is at most $\ucD{D:c:bar} L_k^{-(\alpha-1)/5}$. As for the second term, we have
\begin{align}
    &\PP\left(\mathcal{G},\,|\mathrm{Bar}|\geqslant \frac{M_k}{10},\;|\mathrm{Del}|<\frac{\gamma}{2}\,|\mathrm{Bar}|\right)\nonumber\\
    &=\sum_{M_k/10\leqslant s\leqslant m_k} \PP\left(\mathcal{G},\,|\mathrm{Del}|<\frac{\gamma}{2}s,\,|\mathrm{Bar}|=s\right)\nonumber\\
    &\leqslant \sum_{M_k/10\leqslant s\leqslant m_k} \PP\left(\mathcal{G},\,|\mathrm{Cross}\cap\mathrm{Bar}|>\left(1-\frac{\gamma}{2}\right)s,\,|\mathrm{Bar}|=s\right)&\mbox{using Lemma \ref{D:l:useful}}\nonumber\\
    &\leqslant \sum_{M_k/10\leqslant s\leqslant m_k} \PP\left(\sum_{i=0}^{s-1} \mathbf{1}_{\mathrm{Cross}_{\tau_i}}>\left(1-\frac{\gamma}{2}\right)s\right)\nonumber\\
    &\leqslant \sum_{M_k/10\leqslant s\leqslant m_k} \PP\left(\sum_{i=0}^{s-1}\left( \mathbf{1}_{\mathrm{Cross}_{\tau_i}}-\PP(\mathrm{Cross}_{\tau_i}\vert\mathcal{F}_{\tau_i})\right)>\frac{\gamma}{2}s\right)&\mbox{using \eqref{D:e:lemma2}.}\label{D:e:for_later}
\end{align}
Now, let us set
\begin{align*}
    \left\{\begin{array}{l}
    S_0=0;\\
    \forall l\in\NN^*,\;S_l=\displaystyle\sum_{i=0}^{l-1} R_i,\;\;\text{where}\;\;R_i=\mathbf{1}_{\mathrm{Cross}_{\tau_i}}-\PP(\mathrm{Cross}_{\tau_i}\vert \mathcal{F}_{\tau_i}).
    \end{array}\right.
\end{align*}
Then, $(S_l)_{l\in\NN}$ is a $(\mathcal{F}_{\tau_i})$-martingale: for each $l\in\NN$,
$$\EE[S_{l+1}\vert \mathcal{F}_{\tau_l}]=S_l+\EE[R_l\vert \mathcal{F}_{\tau_l}]=S_l.$$
Also, this martingale has increments bounded by $2$. Therefore we can use Azuma's inequality, which ensures that in the sum in \eqref{D:e:for_later}, we have for each $M_k/10\leqslant s\leqslant m_k$,
\begin{align*}
    \PP\left(\sum_{i=0}^{s-1}\left( \mathbf{1}_{\mathrm{Cross}_{\tau_i}}-\PP(\mathrm{Cross}_{\tau_i}\vert\mathcal{F}_{\tau_i})\right)>\frac{\gamma}{2}s\right)=\PP\left(S_{s}>\frac{\gamma}{2}s\right)\leqslant \exp\left(-\frac{\gamma^2 s}{32}\right)\leqslant e^{-c\,M_k},
\end{align*}
where $c$ does not depend on $y$. Therefore, the sum in \eqref{D:e:for_later} can be bounded from above by $m_k\, e^{-cM_k}\leqslant e^{-cL_k}$, where $c$ does not depend on $y$, hence the result by choosing $\ucD{D:c:del}$ properly.
\end{proof}

We now have all the tools to prove the global delay we are after.
\begin{proposition}
    There exists $\varepsilon>0$ such that
\begin{align}\label{D:e:contradiction'}
p_{\mathbf{H}_k}(v_+-\varepsilon)\xrightarrow[k\to\infty\atop k\in \mathcal{S}]{}0,
\end{align}
\end{proposition}

\begin{proof}
Let us fix $k\in\mathcal{S}$ such that $k>\ucDk{D:k:paths}.$ Assume that $\mathcal{G}_k$ occurs and that for every $y\in I_{\mathbf{H}_k}$, $|\mathrm{Del}(y,k)|\geqslant \frac{\gamma M_k}{20}$. This implies that for every $y\in I_{\mathbf{H}_k}$, for some number $N\geqslant \frac{\gamma M_k}{20r}$ of indexes $j_0\in\llbracket 0,M_k-1\rrbracket$, there exists $j_1\in\llbracket 0,r-1\rrbracket$ such that $\mathrm{Del}_j(y,k)$ occurs, where $j=j_0r+j_1$. For those $N$ indexes, we have a delay at scale $H_k$:
\begin{align*}
    X_{(j_0+1)H_k}^y
    &\leqslant X_{(j+1)h_k}^y+(r-j_1-1)\left(v_++\frac{\delta}{2r}\right)h_k&\mbox{using $\mathcal{G}$ and \eqref{D:e:lam}}\\
    &\leqslant \pi_1(w_j^y)+(v_+-\delta) h_k+(r-j_1-1)\left(v_++\frac{\delta}{2r}\right)h_k&\mbox{since $\mathrm{Del}_j$ occurs}\\
    &\leqslant X^y_{j_0H_k}+j_1v_+h_k+(v_+-\delta) h_k+(r-j_1-1)\left(v_++\frac{\delta}{2r}\right)h_k&\mbox{by \eqref{D:d:w_j}}\\
    &\leqslant X^y_{j_0H_k} +\left(v_+-\frac{\delta}{2r}\right)H_k&\mbox{since $H_k=rh_k$.}
\end{align*}
In the end, using $\mathcal{G}$ in order to bound the displacement for the other indexes $j_0$, we get
\begin{align*}
    X^y_{\mathbf{H}_k}-\pi_1(y)
    &\leqslant N \left(v_+-\frac{\delta}{2r}\right) H_k+ (M_k-N) \left(v_++\lambda\right)H_k\\
    &=M_k(v_++\lambda)H_k-N \left(\frac{\delta}{2r}+\lambda\right)H_k\\
    &\leqslant M_k(v_++\lambda)H_k-\frac{\gamma M_k}{20r} \left(\frac{\delta}{2r}+\lambda\right) H_k\\
    &=\left(v_+-\frac{\gamma\lambda}{20r}\right)\mathbf{H}_k\\
    &<(v_+-\varepsilon)\mathbf{H}_k,
\end{align*} where we set $\varepsilon=\frac{\gamma\lambda}{40r}>0.$ This being true for every $y\in I_{\mathbf{H}_k}$, $A_{\mathbf{H}_k}(v_+-\varepsilon)$ cannot occur. As a result,
\begin{align*}
p_{\mathbf{H}_k}(v_+-\varepsilon)
&\leqslant \PP(\mathcal{G}^c)+\PP\left(\mathcal{G}_k,\,\exists\, y\in I_{\mathbf{H}_k},\,|\mathrm{Del}(y,k)|<\frac{\gamma M_k}{20}\right),\\
&\leqslant \PP(\mathcal{G}^c)+\mathbf{H}_k \sup_{y\in I_{\mathbf{H}_k}} \PP\left(\mathcal{G}_k,\,|\mathrm{Del}(y,k)|<\frac{\gamma M_k}{20}\right)\\
&\leqslant \PP(\mathcal{G}^c)+ L_k^{2}\, \ucD{D:c:del}\,L_k^{-(\alpha-1)/5}&\mbox{using Lemma \ref{D:p:delay}}.
\end{align*}
Now, using \eqref{D:e:proba_G} and the fact that $\alpha>11$, this goes to $0$ when $k$ goes to infinity (with $k\in\mathcal{S}$, where $\mathcal{S}$ is unbounded), which concludes the proof of \eqref{D:e:contradiction'}.
\end{proof}

\begin{proof}[Proof of Lemma \ref{D:l:v_+=v_-}]
Since \eqref{D:e:contradiction'} is in contradiction with the definition of $v_+$, our assumption $v_-<v_+$ is never satisfied, which concludes the proof of Lemma \ref{D:l:v_+=v_-}.
\end{proof}

\section*{Acknowledgments} This work could not have been possible without the extensive help of my PhD supervisors Oriane \textsc{Blondel} (ICJ, Villeurbanne, France) and Augusto \textsc{Teixeira} (IMPA, Rio de Janeiro, Brazil), and I would like to take this opportunity to thank them wholeheartedly for their involvement, kindness and patience. This work was supported by a doctoral contract provided by CNRS. {\color{black}I would like to thank the reviewers who did a very thorough work and provided a great help in improving this article.} My working in person with Augusto \textsc{Teixeira} in IST (Lisbon, Portugal) and IMPA (Rio de Janeiro, Brazil) was enabled by grants from ICJ, Labex Milyon and INSMI.

\bibliographystyle{alpha}
\bibliography{Biblio}

\end{document}